\documentclass[12pt, reqno]{amsart}

\usepackage{amsfonts,amsmath,amsthm,amssymb}
\usepackage{latexsym}

 \DeclareMathOperator{\ord}{ord}
\DeclareMathOperator{\trdeg}{trdeg}
\DeclareMathOperator{\spec}{Spec}

\begin{document}

{\theoremstyle{plain}
  \newtheorem{theorem}{Theorem}[section]
  \newtheorem{corollary}[theorem]{Corollary}
  \newtheorem{proposition}[theorem]{Proposition}
  \newtheorem{lemma}[theorem]{Lemma}
  \newtheorem{question}[theorem]{Question}
  \newtheorem{conjecture}[theorem]{Conjecture}
}

{\theoremstyle{definition}
  \newtheorem{definition}[theorem]{Definition}
  \newtheorem{remark}[theorem]{Remark}
  \newtheorem{example}[theorem]{Example}
}

\numberwithin{equation}{section}

\def\ZZ{\mathbb{Z}}
\def\NN{\mathbb{N}}
\def\R{\mathcal{R}}
\def \a {\alpha}
\def \b {\beta}
\def \d {\delta}
\def \D {\Delta}
\def \e {\epsilon}
\def \l {\lambda}
\def \g {\gamma}
\def \G {\Gamma}
\def \p {\pi}
\def \n {\nu}
\def \m {\mu}
\def \s {\sigma}
\def \S {\Sigma}
\def \w {\omega}
\def \F {\Phi}
\def \O {\mathcal{O}}
\def \W {\Omega}
\def \k {\bold{k}}
\def \U {\mathcal{U}}
\def \t {\tau}
\def \P {\Psi}
\def \et {\eta}
\def \th {\theta}

\def \ra {\rightarrow}

\title[Toroidalization of generating sequences]{Toroidalization of generating
sequences in dimension two function fields}
\author{Laura Ghezzi}
\address{Florida International University, Department of Mathematics,
University Park, Miami, FL 33199, USA} \email{ghezzil@fiu.edu}
\urladdr{http://www.fiu.edu/$\sim$ghezzil/}

\author{Huy T\`ai H\`a}
\address{Tulane University, Department of Mathematics, 6823 St. Charles Ave.,
New Orleans LA 70118, USA} \email{tai@math.tulane.edu}
\urladdr{http://www.math.tulane.edu/$\sim$tai/}

\author{Olga Kashcheyeva}
\address{University of Illinois at Chicago, Department of Mathematics,
Statistics and Computer Science, 851 S. Morgan (m/c 249), Chicago,
IL 60607, USA} \email{olga@math.uic.edu}

\thanks{The authors would like to thank S.D. Cutkosky for suggesting the
problem and for many stimulating discussions on the topic.}
\thanks{The second author is partially supported by the Research Summer
Fellowship of Tulane University.}

\begin{abstract}
Let $\k$ be an algebraically closed field of characteristic 0, and
let $K^*/K$ be a finite extension of algebraic function fields of
transcendence degree 2 over $\k$. Let $\n^*$ be a $\k$-valuation of
$K^*$ with valuation ring $V^*$, and let $\n$ be the restriction of
$\n^*$ to $K$. Suppose that $R \to S$ is an extension of algebraic
regular local rings with quotient fields $K$ and $K^*$ respectively,
such that $V^*$ dominates $S$ and $S$ dominates $R$. We prove that
there exist sequences of quadratic transforms $R \to \bar{R}$ and $S
\to \bar{S}$ along $\n^*$ such that $\bar{S}$ dominates $\bar{R}$
and the map between generating sequences of $\n$ and $\n^*$ has a
toroidal structure. Our result extends the Strong Monomialization
theorem of Cutkosky and Piltant.
\end{abstract}
\maketitle


\section{Introduction}

Let $\k$ be an algebraically closed field of characteristic 0, and
let $K$ be an algebraic function field over $\k$. Throughout this
paper we say that a subring $R$ of $K$ is {\it algebraic} if $R$ is
essentially of finite type over $\k$. We will denote the maximal
ideal of a local ring $R$ by $m_R$.

Let $K^*/K$ be a finite extension of algebraic function fields over
$\k$. Let $\n^*$ be a $\k$-valuation of $K^*$ with valuation ring
$V^*$ and value group $\G^*$. Let $\n$ be the restriction of $\n^*$
to $K$ with valuation ring $V$ and value group $\G$. Consider an
extension of algebraic regular local rings $R \to S$ where $R$ has
quotient field $K$, $S$ has quotient field $K^*$, $R$ is dominated
by $S$ and $S$ is dominated by $V^*$ (i.e., $m_{V} \cap R = m_R$ and
$m_{V^*} \cap S = m_S$).

It has been a classical topic to investigate finite extensions of
rings of algebraic integers and mappings between algebraic curves.
In these cases $K$ and $K^*$ are of transcendence degree 1 over
$\k$, and the homomorphisms of local rings of points are ramified
maps $R \to S$ of discrete (rank 1) valuation rings. We have that $R
= V$ and $S = V^*$ are local Dedekind domains. Suppose that $(u) =
m_R$ and $(x) = m_S$ are the maximal ideals of $R$ and $S$,
respectively, then
\begin{align}
u = x^e \d \label{intro-classical-eq}
\end{align}
where $\d \in S$ is a unit. The corresponding value groups are $\G
\cong \ZZ$ and $\G^* \cong \ZZ$ and we have a natural isomorphism
$\G^*/\G \cong \ZZ_e$.

The study of ramification theory, in general, for valuations in
arbitrary fields was initiated by Krull and further pursued by many
authors (cf. \cite{C&P} and the literature cited there).

In this paper, we are interested in valuations of two dimensional
algebraic function fields, i.e., the situation when $K$ and $K^*$
are of transcendence degree 2 over $\k$. Valuations in dimension two
are completely described by a compact set of data called {\it
generating sequence}. Generating sequences provide a very useful
tool in the study of algebraic surfaces (cf. \cite{Cossart, C&P, GK,
FJ, Spi, T}).

We shall briefly recall the definition of generating sequences, as
in \cite{Spi}. Let ${\G}_+ = \nu(R \backslash \{0\})$ be the
semigroup of $\G$ consisting of the values of nonzero elements of
$R$. For $\gamma \in {\G}_+$, let $I_{\gamma}=\{f\in R \mid \
\n(f)\geq \gamma\}$. A (possibly infinite) sequence $\{ Q_i \}$ of
elements of $R$ is a generating sequence of $\nu$ if for every
$\gamma \in {\G}_+$ the ideal $I_\g$ is generated by the set
$$\{\prod_{i}{Q_i}^{a_i}\mid \ a_i\in \mathbb{N}_0,\
\sum_{i} a_i \n(Q_i)\geq \gamma\}.$$ A generating sequence of $\n$
is {\it minimal} if none of its proper subsequences is a generating
sequence of $\n$. If $\{ Q_i \}$ is a minimal generating sequence
and $\k=V/m_V$ then $\{ \n(Q_i) \} \subset \G$ forms a minimal set
of generators for ${\G}_+$.

  Notice that a
generating sequence $\{ Q_i \}$ in $R$ and the values $\{ \n(Q_i)
\}$ completely describe the valuation $\n$ (see \cite[Definition
1.1]{Spi} and \cite[Section 7.5]{C&P} for more detailed
discussions). A generating sequence of $\n^*$ in $S$ can be defined
similarly.

The aim of this paper is to find structure theorems for generating
sequences of $\n$ and $\n^*$. Our work is inspired by the Strong
Monomialization theorem of Cutkosky and Piltant \cite[Theorem
4.8]{C&P}, which we recall below.

We first need few definitions. Suppose that $R$ is a local domain. A
{\it monoidal} transform $R \to R'$ is a birational extension of
local domains such that $R' = R\big[\frac{I}{x}\big]_m$ where $I$ is
a regular prime ideal of $R$, $0 \not= x \in I$ and $m$ is a prime
ideal of $R\big[\frac{I}{x}\big]$ such that $m \cap R = m_R$. If $I
= m_R$ then the monoidal transform $R \to R'$ is called a {\it
quadratic} transform. In our situation (dimension two) since $R$ is
a regular local ring, any nontrivial monoidal transform $R \to R'$
is a quadratic transform and there exists a regular system of
parameters $(u,v)$ of $R$ such that $R' = R\big[ \frac{u}{v}
\big]_m$, where $m$ is a maximal ideal of $R\big[ \frac{u}{v}
\big]$. We say that $R \to R'$ is a monoidal transform {\it along
$\n$} if $\n$ dominates $R'$.

The celebrated Local Monomialization theorem of Cutkosky
\cite[Theorem 1.1]{Dale2} states that there exist sequences of
monoidal transforms $R \to R_1$ and $S \to S_1$ along $\n^*$ such
that $\n^*$ dominates $S_1$, $S_1$ dominates $R_1$, and there are
regular parameters $(u,v)$ in $R_1$ and $(x,y)$ in $S_1$, units
$\d_1, \d_2 \in S_1$ and a matrix $A = (a_{ij})$ of nonnegative
integers such that $\det A \not= 0$ and
\begin{eqnarray}
\left\{ \begin{array}{rcl} u & = & x^{a_{11}} y^{a_{12}} \d_1 \\
v & = & x^{a_{21}} y^{a_{22}} \d_2. \end{array} \right.
\label{intro-monomialization}
\end{eqnarray}
The existence of $R_1$ and $S_1$ such that
(\ref{intro-monomialization}) holds follows directly from the
standard theorems on resolution
 of singularities, but in general we will not have the essential
  condition that $\det A \not= 0$. The difficulty in Cutkosky's
   work is to achieve the condition $\det A \not= 0$
    (we should note that Cutkosky's Local Monomialization theorem is valid in
arbitrary dimension).

In our situation, under the additional assumption that $\G^*$ is a
non-discrete subgroup of $\mathbb{Q}$ (which is the essential and
subtle case), the Strong Monomialization theorem of Cutkosky and
Piltant \cite[Theorem 4.8]{C&P} further assures that $A$ can be
taken to have the following special form
\begin{eqnarray}
A = \left( \begin{array}{cc} t & 0 \\ 0 & 1 \end{array} \right).
\label{intro-matrix}
\end{eqnarray}
Strong Monomialization is an important and useful result. It shows
that no matter how complicated the structure of the extension $R
\subset S$ is, after blowing-up we obtain a simultaneous resolution,
that is, an extension $R_1 \subset S_1$ of regular local rings such
that $S_1$ is the localization of the integral closure of $R_1$ in
$K^*$, and this extension is very nice and simple, since it is
monomial.

Observe that $u,v \in R_1$ (resp. $x,y \in S_1$) are the first two
members of a generating sequence of $\n$ (resp. $\n^*$). Therefore,
(\ref{intro-matrix}) exhibits a toroidal structure of the first two
elements of such generating sequences.

The definition of toroidal structures of generating sequences of
$\n$ and $\n^*$ is given in Section \ref{statement}. The goal of our
work is to investigate toroidal structures of generating sequences
of $\n$ and $\n^*$. Our main theorem is stated as follows.

\begin{theorem}[Theorem \ref{monomialization}] \label{intro-main}
Let $\k$ be an algebraically closed field of characteristic 0, and
let $K^*/K$ be a finite extension of algebraic function fields of
transcendence degree 2 over $\k$. Let $\n^*$ be a $\k$-valuation of
$K^*$ with valuation ring $V^*$, and let $\n$ be the restriction of
$\n^*$ to $K$. Suppose that $R \to S$ is an extension of algebraic
regular local rings with quotient fields $K$ and $K^*$ respectively,
such that $V^*$ dominates $S$ and $S$ dominates $R$. Then there
exist sequences of quadratic transforms $R \to \bar{R}$ and $S \to
\bar{S}$ along $\n^*$ such that $\bar{S}$ dominates $\bar{R}$ and
the map between generating sequences of $\n$ and $\n^*$ in $\bar{R}$
and $\bar{S}$ respectively, has a toroidal structure.
\end{theorem}

To prove this theorem, we consider different cases according to
Zariski's classification of valuations in two dimensional function
fields over an algebraically closed field of characteristic 0
\cite{Z}. In most cases, the result follows from a standard
application of the Strong Monomialization theorem of Cutkosky and
Piltant. These cases are analyzed in Section \ref{easycases}. The
bulk of the paper is devoted to the essential case, when $\G^*$ is a
non-discrete subgroup of $\mathbb{Q}$. We shall now briefly describe
the main steps of the proof in this case.

Let $(x,y)$ be a regular system of parameters in $S$. We begin in
Section \ref{jumpingpoly} by constructing a sequence of {\it jumping
polynomials} $\{ T_i \}_{i \geq 0}$ in $S$ with $T_0 = x$ and $T_1 =
y$, which contains as a subsequence a minimal generating sequence of
$\n^*$. Our notion of jumping polynomials is very similar to Favre
and Jonsson's notion of {\it key polynomials} \cite{FJ}. The idea of
key polynomials is originally due to MacLane \cite{MacL}. By
normalizing, we may assume that $\n^*(x) = 1$. Let $\n^*(y) =
\frac{p_1}{q_1}$, where $p_1$ and $q_1$ are coprime positive
integers. For each $i \ge 1$, we define $T_{i+1}$ recursively. Let
$p_{i+1}$ and $q_{i+1}$ be coprime positive integers defined by
$$\n^*(T_{i+1}) = q_i \n^*(T_i) + \dfrac{1}{q_1 \dots q_i}\cdot
\dfrac{p_{i+1}}{q_{i+1}}.$$

\noindent Our proof of Theorem \ref{intro-main} proceeds in the
following line of arguments.

\begin{enumerate}
\item[(1)] We observe that the collection of jumping polynomials
$\{ T_i \}_{i\geq 0}$ forms a generating sequence of $\n^*$ in $S$
(Theorem \ref{nondiscretesequence}).

\item[(2)] Let $d = d(R,S)$ be the greatest common divisor of $\{
p_i ~|~ i \ge 1 \}$. We show that the powers of $x$ that appear in
$T_i$ are multiple of $d$ for all $i \ge 2$ (Corollary
\ref{xpowers}). In other words $T_i$, for $i \ge 2$, is a polynomial
in $x^d$ and $y$.

\item[(3)] Let us denote by $t(R,S)$ the power $t$ in
(\ref{intro-matrix}) obtained from the Strong Monomialization
theorem. For simplicity, assume that the constant $\d_1$ of
(\ref{intro-monomialization}) is equal to 1. We observe that if
$t(R,S)$ divides $d$ then $T_i$, for $i \ge 2$, is a polynomial in
$u$ and $v$. This shows that $u$, together with the collection
$\{T_i\}_{i\geq 1}$, form a generating sequence of $\n$ in $R$.
Therefore we obtain the desired toroidal structure.

\item[(4)] The core of our argument is to show that if $t(R,S)$
does not divides $d(R,S)$ then we can find sequences of quadratic
transforms $R \to R'$ and $S \to S'$ so that $t(R',S') < t(R,S)$.
More precisely, let $M=\min\{i>0|\;t\nmid p_i\}$. Then $t(R',S')$ is
the greatest common divisor of $t$ and $p_M$. Lemma \ref{chunk} is
crucial in the proof of this step.

Finally, starting with a sequence of jumping polynomials in $S'$ and
repeating the above process, after a finite number of iterations we
end up with the situation where $t(R',S')$ divides $d(R',S')$. Then
we conclude as in (3).

\end{enumerate}

We remark that in order to make our arguments work we need a very
explicit description of the quadratic transforms that we perform.
There are several preparatory lemmas to this end.


\section{Statement of the result} \label{statement}

Let $\k$ be an algebraically closed field of characteristic 0 and
let $K^*/K$ be a finite extension of algebraic function fields of
transcendence degree 2 over $\k$. Let $\n^*$ be a $\k$-valuation of
$K^*$ with valuation ring $V^*$ and value group $\G^*$ and let $\n$
be the restriction of $\n^*$ to $K$ with valuation ring $V$ and
value group $\G$.

Suppose that $S$ is an algebraic regular local ring with quotient
field $K^*$ which is dominated by $V^*$ and $R$ is an algebraic
regular local ring with quotient field $K$ which is dominated by
$S$. We will show that there exist sequences of quadratic transforms
$R\ra R'$ and $S\ra S'$ along $\n^*$ such that $S'$ dominates $R'$
and the map between generating sequences of $S'$ and $R'$ has the
following toroidal structure.

\begin{itemize}

\item[(1)] If $\n^*$ is divisorial then $R'=V$ and $S'=V^*$ with
regular parameters $u\in R'$ and $x\in S'$ such that $u=x^a\g$ for
some unit $\g\in S'$. We also have that $\{x\}$ is a minimal
generating sequence of $\n^*$ and $\{u\}$ is a minimal generating
sequence of $\n$.

\item[(2)] If $\n^*$ has rank 2 then there exist regular
parameters $(x,y)$ in $S'$ and $(u,v)$ in $R'$ such that $\{x,y\}$
is a minimal generating sequence of $\n^*$, $\{u,v\}$ is a minimal
generating sequence of $\n$ and
\begin{align*}
u &=x^ay^b\d\\
v &=y^d\g
\end{align*}
for some units $\d,\g\in S'$, and for some nonnegative integers $a,
b, d$ such that $ad\neq 0$.

\item[(3)] If $\n^*$ has rank 1 and rational rank 2 then there
exist regular parameters $(x,y)$ in $S'$ and $(u,v)$ in $R'$ such
that $\{x,y\}$ is a minimal generating sequence of $\n^*$, $\{u,v\}$
is a minimal generating sequence of $\n$ and
\begin{align*}
u &=x^ay^b\d\\
v &=x^cy^d\g
\end{align*}
for some units $\d,\g\in S'$, and for some nonnegative integers $a,
b, c, d$ such that $ad-bc \neq 0$.

\item[(4)] If $\G$ and $\G^*$ are non-discrete subgroups of
$\mathbb{Q}$ then there exist a minimal generating sequence
$\{H_i\}_{i\ge 0}$ of $\n^*$ in $S'$ and regular parameters $(u,v)$
in $R'$ such that
\begin{align*}
u &=H_0^a\g\\
v &=H_1
\end{align*}
for some unit $\g\in S'$, and $H_i\in R'$ for all $i>1$.
Furthermore, $\{u,\{H_i\}_{i>0}\}$ is a generating sequence of $\n$
in $R'$.

\item[(5)] If $\n$ is discrete but not divisorial then there exist
regular parameters $(x,y)$ in $S'$ and $(u,v)$ in $R'$ such that
$\G^*$ is generated by $\n^*(x)$, $\G$ is generated by $\n(u)$ and
$u=x^a\g$ for some unit $\g\in S'$. Moreover, $S'$ has a non-minimal
generating sequence $\{x,\{T_i\}_{i>0}\}$  such that
$\{u,\{T_i\}_{i>0}\}$ form a non-minimal generating sequence in
$R'$.

\end{itemize}


\section{Valuations in 2 dimensional function fields} \label{easycases}

Zariski in \cite{Z} gave a classification of valuations in two
dimensional function fields over an algebraically closed field of
characteristic zero. We refer to \cite{Dalebook} (Chapter 8, Section
1) for a modern treatment of the subject and for the definitions and
background needed in this section.

We will prove our main theorem by analyzing the different types of
valuations of $K^*$. Notations are as in Section \ref{statement}.

\subsection{One dimensional valuations}

By definition, $\n^*$ is divisorial. In this case $\n$ and $\n^*$
are discrete, and $V$ and $V^*$ are iterated quadratic transforms of
$R$ and $S$ respectively (see \cite[Proposition 4.4]{Ab}).

Let $u$ be a regular parameter of $V$ and let $x$ be a regular
parameter of $V^*$. Then there is a relation $$u=x^a\g$$ where
$\g\in V^*$ is a unit, and $a\geq 1$. Since $\{ u \}$ is a mimimal
generating sequence for $V$, and $\{ x \}$ is a minimal generating
sequence for $V^*$ the theorem is proved.

\subsection{Zero dimensional valuations of rational rank 2}

By Strong Monomialization \cite[Theorem 4.8]{C&P} there exist
sequences of quadratic transforms $R\ra R'$ and $S\ra S'$ along
$\n^*$ such that $R'$ has regular parameters $(u,v)$, $S'$ has
regular parameters $(x,y)$, and
\begin{align*}
u &=x^ay^b\d\\
v &=x^cy^d\g
\end{align*}
for some units $\d,\g\in S'$ and for some nonnegative integers $a,
b, c, d$ such that $ad-bc \neq 0$. Further, $c=0$ if $\n^*$ has rank
two. We also have that $\{\n^*(x), \n^*(y)\}$ is a rational basis of
$\G^*\otimes \mathbb{Q}$, and $\{\n(u), \n(v)\}$ is a rational basis
of $\G\otimes \mathbb{Q}$.

Let $z\in S'$. Then $z\in \hat S'=S'/m_{S'}[[x,y]]=\k[[x,y]]$, since
$\n^*$ is zero dimensional and $\k$ is algebraically closed. Observe
that $z$ has an expansion $z=\sum_{i\geq 1}a_ix^{b_i}y^{c_i}$, where
$a_i\in \k$, $b_i$ and $c_i$ are non negative integers, and the
terms have increasing value, since $\n^*(x)$ and $\n^*(y)$ are
rationally independent. It follows that
$\n^*(z)=b_1\n^*(x)+c_1\n^*(y)$. Hence $\{ x,y \}$ is a minimal
generating sequence of $\n^*$ in $S'$, and similarly $\{ u,v \}$ is
a minimal generating sequence of $\n$ in $R'$, and the theorem is
proved.

\medskip

The rest of the paper will be devoted to studying the remaining
cases, that is zero dimensional valuations of rational rank 1.

\subsection{Non-discrete zero dimensional valuations of rational rank 1}

We can normalize $\G^*$ so that it is an ordered subgroup of
$\mathbb{Q}$, whose denominators are not bounded, as $\G^*$ is not
discrete. In Example 3, Section 15, Chapter VI of \cite{ZS},
examples are given of two-dimensional algebraic function fields with
value group equal to any given subgroup of the rational numbers.
This case is much more subtle.

\subsection{Discrete zero dimensional valuations of rational rank 1}

If $\n^*$ is discrete, then $\n$ is also discrete. This case will be
handled in the same way as the case of non-discrete zero dimensional
valuations of rational rank 1, but the generating sequences of
$\n^*$ and $\n$ will not be minimal.



\section{Construction of jumping polynomials}\label{jumpingpoly}

From now on we work under the assumption that the value group of
$\n^*$ is a subgroup of $\mathbb{Q}$ and
$\trdeg_{\k}(V^*/{m_{V^*}})=0$. Let $(x,y)$ be a system of regular
parameters in $S$. We normalize the value group $\G^*$ of $K^*$ so
that $\n^*(x)=1$.

We shall now construct a sequence of polynomials $\{T_i\}_{i\geq 0}$
in $S$. Let
\begin{align*}
\left\{ \begin{array}{ccc} T_0 & = & x\\
T_1 & = & y. \end{array} \right.
\end{align*}
Set $q_0=\infty$ and choose a pair of coprime positive integers
$(p_1,q_1)$ so that $\n^*(y)=\dfrac{p_1}{q_1}$. For $i \ge 1$,
$T_{i+1}$ is defined recursively as follows.
 Let
$$T_{i+1}=T_i^{q_i}-\l_i\prod_{j=0}^{i-1}T_j^{n_{i,j}},$$
where $n_{i,j}<q_j$ is a nonnegative integer such that
$q_i\nu^*(T_i)=\nu^*(\prod_{j=0}^{i-1}T_j^{n_{i,j}})$, that is
$q_i\nu^*(T_i)=\sum_{j=0}^{i-1}n_{i,j}\n^*(T_j)$, and $\l_i\in \k$
is the residue of
$\dfrac{T_i^{q_i}}{\prod_{j=0}^{i-1}T_j^{n_{i,j}}}$.

Finally, choose positive integers $p_{i+1}$ and $q_{i+1}$ so that
$(p_{i+1},q_{i+1})=1$ and
$$\nu^*(T_{i+1})=q_i\nu^*(T_i)+\dfrac{1}{q_1\cdots
q_i}\cdot\dfrac{p_{i+1}}{q_{i+1}}.$$

\begin{definition}\label{jp} The polynomial $T_{i}$ will be called the $i$-th
{\it jumping polynomial} and the value $\n^*(T_i)$ will be called
the $i$-th j-{\it value}. We will denote the $i$-th j-value by
$\b_i$. We say that $\b_i$ is an {\it independent} j-value if
$q_i\neq 1$.

\end{definition}

\begin{remark}\label{b-inequality}
For $i>0$ let $Q_i=q_1\cdots q_i$. Observe that $Q_i\b_i$ is an
integer number,
$\b_{i+1}=q_i\b_i+\dfrac{1}{Q_i}\cdot\dfrac{p_{i+1}}{q_{i+1}}$ and
$q_{i+1}\b_{i+1}\geq\b_{i+1}>q_i\b_i$.
\end{remark}

Consider the subsequence $\{\b_{i_l}\}_{l\geq 0}$ of all independent
j-values. Let $\bar{\b}_l=\b_{i_l}$ denote the $l$-th independent
j-value, $\bar{q}_l=q_{i_l}$ and $H_l=T_{i_l}$. Since
$$
T_{i+1}=T_i^{q_i}-\l_i\prod_{j=0}^{i-1}T_j^{n_{i,j}},\text { where }
0\leq n_{i,j}<q_j,
$$
it follows that $n_{i,j}=0$ whenever $q_j=1$. Therefore only the
$H_l$'s with $i_l < i$ will appear in the product
$\prod_{j=0}^{i-1}T_j^{n_{i,j}}$. Thus, if $0<i=i_l$ then
$$
T_{i+1}=H_l^{\bar{q}_l}-\l_i\prod_{j=0}^{l-1}H_j^{n_{i,i_j}}.
$$
If $i+1 <i_{l+1}$ then $q_{i+1}=1$ and
$$
T_{i+2}=T_{i+1}-\l_{i+1}\prod_{j=0}^l H_j^{n_{i+1,i_j}}=
H_l^{\bar{q}_l}-\l_i\prod_{j=0}^{l-1}H_j^{n_{i,i_j}}-
\l_{i+1}\prod_{j=0}^l H_j^{n_{i+1,i_j}}.
$$
In general, the recursive formula for $H_{l+1}$ with $l>0$ will be
\begin{align*}
H_{l+1} &
=H_l^{\bar{q}_l}-\l_{i_l}\prod_{j=0}^{l-1}H_j^{n_{i_l,i_j}}-
\l_{i_l+1}\prod_{j=0}^l H_j^{n_{i_l+1,i_j}}-\l_{i_l+2}\prod_{j=0}^l
H_j^{n_{i_l+2,i_j}}- \dots\\
& \quad \dots -\l_{i_{l+1}-1}\prod_{j=0}^l H_j^{n_{i_{l+1}-1,i_j}}
=H_l^{\bar{q}_l}-\l_{i_l}\prod_{j=0}^{l-1}H_j^{n_{i_l,i_j}}-
\sum_{i'=i_l+1}^{i_{l+1}-1}\l_{i'}\prod_{j=0}^l H_j^{n_{i',i_j}}.
\end{align*}

We also notice that the sequence of independent jumping polynomials
$\{H_l\}_{l\geq 0}$ starts with $H_0=x$ and
$H_1=y-\sum_{j=1}^{i_1-1}\l_jx^{\b_j}$.

Independent j-values furthermore have a number of basic properties.
If $l>0$ and $i_l\le i<i_{l+1}$ then $q_1\cdots
q_i=\bar{q}_1\cdots\bar{q}_l$ and the following equalities hold
\begin{align*}
\bar{\b}_0 &=\b_0=1,\;\bar{q}_0=q_0=\infty,\\
\bar{\b}_1 &=p_1+\dots+p_{i_1-1}+\dfrac{p_{i_1}}{q_{i_1}},\\
\bar{\b}_{l+1} &=\bar{q}_l\bar{\b}_l+
\dfrac{p_{i_l+1}+p_{i_l+2}+\dots+p_{i_{l+1}-1}}{\bar{q}_1\cdots\bar{q}_l}+
\dfrac{1}{\bar{q}_1\cdots\bar{q}_l}\cdot\dfrac{p_{i_{l+1}}}{\bar{q}_{l+1}}.
\end{align*}

\begin{remark}\label{b*-inequality}
For all $l>0$ denote by $\bar{Q}_l=\bar{q}_1\cdots\bar{q}_l$ and
$\bar{p}_l=(p_{i_{l-1}+1}+\dots+p_{i_l-1})\bar{q}_l+p_{i_l}$. Then
$(\bar{p}_l,\bar{q}_l)=(p_{i_l},q_{i_l})=1$,
$\bar{\b}_1=\dfrac{\bar{p}_1}{\bar{q}_1}$ and
$\bar{\b}_{l+1}=\bar{q}_l\bar{\b}_l+\dfrac{1}{\bar{Q}_l}
\cdot\dfrac{\bar{p}_{l+1}}{\bar{q}_{l+1}}$. In particular,
$\bar{q}_{l+1}\bar{\b}_{l+1}>\bar{\b}_{l+1}>\bar{q}_l\bar{\b}_l$.
\end{remark}

\begin{remark}\label{normalize} In general, if $(x,y)$ is a system of regular
parameters in $S$ we may not necessarily have $\n^*(x)=1$. Then in
order to define a sequence of jumping polynomials $\{T_i\}_{i\ge 0}$
such that $T_0=x$ and $T_1=y$, we introduce the following valuation
$\tilde\n$ of $K^*$
$$
\tilde\n(f)=\frac{\n^*(f)}{\n^*(x)}
$$
for all $f\in K^*$. Then $\tilde\n(x)=1$ and we use the construction
above with $\n^*$ replaced by $\tilde\n$. This procedure is
equivalent to normalizing the value group $\G^*$ so that
$\n^*(x)=1$.
\end{remark}

We will see in Section \ref{arithmetic} that the sequence of jumping
polynomials $\{T_i\}_{i\ge 0}$ in $S$ is well defined (Corollary
\ref{welldefined}). The next goal is to show that it forms a
generating sequence of $\n^*$.

\subsection{Discrete case.}
We suppose that the value group of $\n^*$ is isomorphic to $\ZZ$.
After performing a sequence of quadratic transforms along $\n^*$ and
normalizing $\G^*$ we may suppose that $S$ has a system of regular
parameters $(x,y)$ such that $\n^*(x)=1$ generates $\G^*$. Then by
\cite{Spi} (p.154) we have that a set $\{Q_i\}_{i\ge 0}\subset S$ is
a generating sequence of $\n^*$ provided $\n^*(Q_0)=1$, each $Q_i$
is a regular parameter of $S$ such that $(Q_0,Q_i)$ form a system of
regular parameters, and $\lim_{i\ra\infty}\n^*(Q_i)=\infty$.

In particular, there are no minimal generating sequences in $S$. Any
infinite subsequence of a generating sequence which contains $Q_0$
is a generating sequence itself.

\begin{theorem}\label{discretesequence}
The above $\{T_i\}_{i\ge 0}\subset S$ form a generating sequence of
$\n^*$.
\end{theorem}

\begin{proof}
Since $\n^*(x)$ generates $\G^*$ we see that $q_i=1$ for all $i\ge
1$. Thus $T_1=y$ and
$T_{i+1}=y-\l_1x^{n_{1,0}}-\l_2x^{n_{2,0}}-\dots-\l_ix^{n_{i,0}}$
are linear in $y$ for all $i\ge 1$. In particular $T_i$ is a regular
parameter of $S$ and $(x,T_i)$ form a system of regular parameters
in $S$. Notice also that $\b_i=\n^*(T_i)\in\ZZ$ and
$\b_{i+1}>\b_{i}$. This implies that
$\lim_{i\ra\infty}\n^*(T_i)=\infty$.
\end{proof}

\subsection{Non-discrete case.}
We assume now that the value group of $\n^*$ is a non-discrete
subgroup of $\mathbb{Q}$.

\begin{theorem}\label{nondiscretesequence}
With notations as above, $\{T_i\}_{i\ge 0}\subset S$ form a
generating sequence of $\n^*$. Furthermore, if $p_1>1$ then
$\{H_l\}_{l\ge 0}\subset S$ form a minimal generating sequence of
$\n^*$. If $p_1=1$ then $\{H_l\}_{l\ge 1}\subset S$ form a minimal
generating sequence of $\n^*$.
\end{theorem}

\begin{proof}
This is shown in \cite{Spi}, in Chapter 2 of \cite{FJ}, and in
\cite{Cossart}.
\end{proof}

We will give an alternative proof of the above theorem in Section
\ref{remgenseq}. Namely, we will show that the set $\{H_l\}_{l\ge
0}\subset S$ satisfies the sufficient condition for a sequence of
elements of $S$ to be a generating sequence of $\n^*$ given in
\cite{C&P}.

\begin{remark}\label{genvalues}
The two preceding theorems imply that the values of jumping
polynomials $\{\n^*(T_i)\}_{i\ge 0}$ generate the value group
$\G^*$.
\end{remark}

\section{Arithmetics} \label{arithmetic}

In this section we prove several properties of the numbers $q_i$ and
$\b_i$ defined in Section \ref{jumpingpoly}.

\begin{definition}
Given two rational numbers $a$ and $b$ we say that $a$ is {\it
$\mathbb{Z}$-divisible} by $b$, or equivalently $b$ {\it
$\ZZ$-divides} $a$, and write $b | a$, if $a$ is an integer multiple
of $b$, that is $a\in b\mathbb{Z}$ or $a=nb$ for some
$n\in\mathbb{Z}$. The {\it greatest common divisor} of $a$ and $b$,
denoted by $(a,b)$, is as usual the greatest rational number $g$
such that $g|a$ and $g|b$.
\end{definition}


\begin{proposition}
Let $p,q,t$ be nonzero integers with $(p,q)=1$. Then
$(\frac{1}{t},\frac{p}{tq})=\frac{1}{tq}$.
\end{proposition}
\begin{proof}
If $g=(\frac{1}{t},\frac{p}{tq})$ then $\frac{1}{tq}|g$. On the
other hand
$$
\frac{1}{tq}=\frac{1}{tq}(\a p+\b q)=\a\frac{p}{tq}+\b\frac{1}{t}
$$
for some integers $\a$ and $\b$ since $1=(p,q)$. Thus $g
\big|\frac{1}{tq}$.
\end{proof}

\begin{proposition} For $k \ge 1$, we have $(\b_0,\b_1,\dots,\b_k) =
\frac{1}{Q_k}$.
\end{proposition}

\begin{proof}
We use induction on $k$. For $k = 1$, clearly $(\b_0,\b_1) =
(1,\frac{p_1}{q_1}) = \frac{1}{q_1}$. Assume now that
$(\b_0,\b_1,\dots,\b_{k-1})= \frac{1}{Q_{k-1}}$. Then
$$(\b_0,\b_1,\dots,\b_{k})=((\b_0,\b_1,\dots,\b_{k-1}),\b_{k}-q_{k-1}\b_{k-1})
=(\frac{1}{Q_{k-1}},\frac{p_k}{Q_{k-1}q_k}) =\frac{1}{Q_k}.$$
\end{proof}

\begin{corollary}\label{groupvalue}
For $k \ge 0$, let $\G_k = \langle \b_0,\b_1,\dots,\b_k \rangle$.
Then, $\G_k = \frac{1}{Q_k}\mathbb{Z}$ for all $k \ge 1$. {\rm(}That
is, the group generated by the values of the first $k+1$ jumping
polynomials is isomorphic to $\frac{1}{Q_k}\mathbb{Z}$.{\rm)}
\end{corollary}

\begin{proposition} For $k \ge 1$, we have
$q_k\b_k\in\G_{k-1}$. Moreover, if $q_k>1$ then $q_k\b_k$ has order
$q_k$ in $\frac{\G_{k-1}}{q_k\G_{k-1}}$.
\end{proposition}

\begin{proof} We have
$q_k\b_k=q_kq_{k-1}\b_{k-1}+p_k\frac{1}{Q_{k-1}}$ is
$\mathbb{Z}$-divisible by $\frac{1}{Q_{k-1}}$. Thus,
$q_k\b_k\in\G_{k-1}$.

Moreover, if $q_k>1$ then
$G=\frac{\G_{k-1}}{q_k\G_{k-1}}\cong\mathbb{Z}_{q_k}$ is not trivial
and $\ord_G(q_k\b_k)=\ord_G(p_k \frac{1}{Q_{k-1}})
=\ord_{\mathbb{Z}_{q_k}}(p_k)=q_k$ since $(p_k,q_k)=1$.
\end{proof}

\begin{corollary}\label{irreducibility requirement}
For $k \ge 1$, we have
$(\bar{\b}_0,\bar{\b}_1,\dots,\bar{\b}_k)=\frac{1}{\bar{Q}_k}$ and
$\bar{\G}_k= \langle \bar{\b}_0,\bar{\b}_1,\dots,\bar{\b}_k \rangle
= \frac{1}{\bar{Q}_k}\mathbb{Z}$. Also,
$\bar{q}_k\bar{\b}_k\in\bar{\G}_{k-1}$ and $\bar{q}_k\bar{\b}_k$ has
order $\bar{q}_k$ in
$\frac{\bar{\G}_{k-1}}{\bar{q}_k\bar{\G}_{k-1}}$.
\end{corollary}

\begin{remark} With notations as above, we have
$\bar{\G}_k=\G_{i_k}$ and $\bigcup_{k\ge 0}\bar{\G}_k=\bigcup_{k\ge
0}\G_k$.
\end{remark}

\begin{proposition}\label{unique representation}
If $x\in\G_k$ and $x\geq q_k\b_k$ then there exists a unique
representation
\begin{equation}\label{*}
 x=\sum_{j=0}^ka_j\b_j
\end{equation}
 with integer coefficients $0\leq a_j<q_j$.
\end{proposition}

\begin{proof}
We first show existence of the presentation (\ref{*}). We use
induction on $k$. The claim is trivial for $k=0$. Let $x\in\G_1$ and
$x\geq p_1$, that is $x=y+\frac{p}{q_1}$ for some $p,y\in\mathbb{Z}$
such that $y\geq p_1$ and $0\leq p<q_1$. Since $(p_1,q_1)=1$ there
exists an integer $0\leq a_1<q_1$ such that $a_1p_1=p+tq_1$ for some
$t\in\mathbb{Z}$. Notice that $a_1p_1<q_1p_1$ and therefore $t<p_1$.
So $x=(y-t)\b_0+a_1\b_1$ is the required presentation.

Now assume that $k\geq 2$ and that a presentation (\ref{*}) exists
for $k-1$. Let $x\in\G_k$ and $x\geq q_k\b_k$, then
$x=\frac{y}{Q_{k-1}}+\frac{p}{Q_{k-1}q_k}$ for some
$p,y\in\mathbb{Z}$ such that $0\leq p<q_k$. Since $(p_k,q_k)=1$
there exists an integer $0\leq a_k<q_k$ such that $a_kp_k=p+tq_k$
for some $t\in\mathbb{Z}$. Then
$x=\frac{y}{Q_{k-1}}+a_k(\b_k-q_{k-1}\b_{k-1})-\frac{t}{Q_{k-1}}$,
so that
$x-a_k\b_k=\frac{y-t}{Q_{k-1}}-a_kq_{k-1}\b_{k-1}\in\G_{k-1}$ and
$x-a_k\b_k\geq (q_k-a_k)\b_k\geq\b_k>q_{k-1}\b_{k-1}$. Thus, by the
inductive assumption we have that
$x-a_k\b_k=\sum_{j=0}^{k-1}a_j\b_j$ with integer coefficients $0\leq
a_j<q_j$.

To prove uniqueness of the presentation (\ref{*}) it suffices to
show that if $\sum_{j=0}^k c_j\b_j=0$ for some integer coefficients
$-q_j<c_j<q_j$ then $c_j=0$ for all $j$. We again use induction on
$k$. The claim is trivial for $k=0$. Assume that the claim is true
for $k-1$ and suppose that $\sum_{j=0}^k c_j\b_j=0$ for some integer
coefficients $-q_j<c_j<q_j$. We only need to show that $c_k=0$.
Since $c_k\b_k=-\sum_{j=0}^{k-1}c_j\b_j$, $c_k\b_k\in\G_{k-1}$ and
therefore $c_k\frac{p_k}{Q_{k-1}q_k}\in\G_{k-1}$. Thus
$\frac{c_kp_k}{Q_{k-1}q_k}$ is $\mathbb{Z}$-divisible by
$\frac{1}{Q_{k-1}}$, that is $q_k|c_kp_k$. This implies that
$q_k|c_k$, since $(p_k,q_k)=1$, and therefore $c_k=0$ due to the
inequality $-q_k<c_k<q_k$.
\end{proof}

\begin{corollary}\label{justification}
For all $k>0$ there exists a unique representation
$q_k\b_k=\sum_{j=0}^{k-1}n_{k,j}\b_j$ with integer coefficients
$0\leq n_{k,j}<q_j$.
\end{corollary}

\begin{proof}
The statement is clear for $k=1$, since $q_1\b_1=p_1\b_0$. For $k>1$
the conclusion follows from Proposition \ref{unique representation},
since $q_k\b_k\in \G_{k-1}$ and $q_k\b_k > q_{k-1}\b_{k-1}$.
\end{proof}

\begin{corollary}\label{welldefined}
In the notations of Section \ref{jumpingpoly}, the sequence of
jumping polynomials $\{T_i\}_{i\geq 0}$ in $S$ is well-defined.
\end{corollary}

\begin{proof}
The statement follows immediately from Corollary \ref{justification}
and the assumption that $V^*/{m_{V^*}}=\k$.
\end{proof}





We now recall some well-known facts about continued fractions. Let
$p$ and $q$ be positive integers such that $(p,q)=1$. Consider the
Euclidian algorithm for finding the greatest common divisor of $p$
and $q$:
\begin{align*}
r_0 &=f_1r_1+r_2\\
r_1 &=f_2r_2+r_3\\
&\dots\\
r_{N-2} &=f_{N-1}r_{N-1}+1\\
r_{N-1} &=f_N\cdot 1,
\end{align*}
where $r_0=p$, $r_1=q$ and $r_1>r_2>\dots>r_{N-1}>r_N=1$. Denote by
$N=N(p,q)$ the number of divisions in the Euclidian algorithm for
$p$ and $q$ and by $f_1,\,f_2,\dots,f_N$ the coefficients in the
Euclidian algorithm for $p$ and $q$. Define $F_i=f_1+\dots+f_i$ and
$\e(p,q)=f_1+\dots+f_N=F_N$,
$f_1(p,q)=f_1=\left[\frac{p}{q}\right]$. Let $a$ and $b$ be integers
such that $0<a\le p$, $0\le b<q$, and $aq-bp=1$.

\begin{remark} With notations as above,
$ \dfrac{p}{q}=f_1+\dfrac{1}{f_2+\dots+\dfrac{1}{f_N}}.$
\end{remark}

Let $\{P_k(z_1,\dots,z_k)\}_{k\in \mathbb{N}_0}$ be a sequence of
polynomials as in \cite{Spi}. So
$P_k(z_1,\dots,z_k)\in\mathbb{N}_0[z_1,\dots,z_k]$ is a polynomial
in $k$ variables with nonnegative integer coefficients such that for
any set of numbers ${c_1,\dots,c_K}$ we have
$$
c_1+\dfrac{1}{c_2+\dots
+\dfrac{1}{c_K}}=\dfrac{P_K(c_1,\dots,c_K)}{P_{K-1}(c_2,\dots,c_K)}.
$$
We also assume that $P_0=1$ and set $P_{-1}=0$.

Then it follows from properties (1.2)-(1.6) in \cite{Spi} that
\begin{align*}
p&=P_N(f_1,\dots,f_N),\\
q&=P_{N-1}(f_2,\dots,f_N),\\
a&=P_{N-1}(f_1,\dots,f_{N-1}),\quad\quad\quad
b=P_{N-2}(f_2,\dots,f_{N-1}),
\quad\;\quad \text{ if $N$ is odd,}\\
a&=p-P_{N-1}(f_1,\dots,f_{N-1}),\quad
b=q-P_{N-2}(f_2,\dots,f_{N-1}),\quad \text{ if $N$ is even.}
\end{align*}

We also recall property (1.5) from \cite{Spi} here since it will be
used in the sequel
\begin{align*}
P_k(f_1,\dots,f_k)
&=f_kP_{k-1}(f_1,\dots,f_{k-1})+P_{k-2}(f_1,\dots,f_{k-2}),\\
P_{k-1}(f_2,\dots,f_k)
&=f_kP_{k-2}(f_2,\dots,f_{k-1})+P_{k-3}(f_2,\dots,f_{k-2}).
\end{align*}


\section{Sequences of quadratic transforms} \label{quadratic transforms}

We will now consider a sequence
$$
S=S_0\ra S_1\ra S_2\ra\dots\ra S_i\ra\dots
$$
of quadratic transforms along $\n^*$. Suppose that $E$ is a
nonsingular irreducible curve on $\spec S$. Denote by $\p_i$ the map
$\spec S_i\ra\spec S$ and by $E_i$ the reduced simple normal
crossing divisor $\p_i^{-1}(E)_{red}$. We say that $S_i$ is {\it
free} if $E_i$ has exactly one irreducible component. For a free
ring $S_i$ and a regular parameter $x_i\in S_i$ we will say that
$x_i$ is an {\it exceptional parameter} if $x_i$ is supported on
$E_i$. A system of parameters $(x_i,y_i)$ of a free ring $S_i$ is
called {\it permissible} if $x_i$ is an exceptional parameter.

If $S_i$ has regular parameters $(x_i,y_i)$ then we can choose
regular parameters $(x_{i+1},y_{i+1})$ in $S_{i+1}$ as follows
\begin{itemize}
\item[a)] if $\n^*(x_i)<\n^*(y_i)$ then
$ x_{i+1}=x_i$ and $y_{i+1}=\frac{y_i}{x_i}$,
\item[b)] if $\n^*(x_i)>\n^*(y_i)$ then
$ x_{i+1}=\frac{x_i}{y_i}$ and $y_{i+1}=y_i$,
\item[c)] if $\n^*(x_i)=\n^*(y_i)$ then
$ x_{i+1}=x_i$ and $y_{i+1}=\frac{y_i}{x_i}-c$, where $c\in\k$ is
the residue of $\frac{y_i}{x_i}$.
\end{itemize}

Our goal is to describe explicitly the sequence of quadratic
transforms of $S$ along $\n^*$. Assume that $(x,y)$ is a permissible
system of parameters in $S$. Let $p$ and $q$ be positive coprime
integers such that $\frac{\n^*(y)}{\n^*(x)}=\frac{p}{q}$. We denote
by $\m$ the value $\n^*(x)$. Let $N=N(p,q)$, $f_1,\dots,f_n$ and
$F_1,\dots, F_N$ be defined by the Euclidian algorithm for $p$ and
$q$ as in Section \ref{arithmetic}. Let $a$ and $b$ be integers such
that $0<a\le p$, $0\le b<q$ and $aq-bp=1$. We will investigate the
following sequence of quadratic transforms along $\n^*$
$$
S=S_0\ra S_1\ra\dots\ra S_{F_1}\ra\dots\ra S_{F_j}\ra\dots\ra
S_{F_N}.
$$
If $N>1$ then for all $0\le j\le F_1$, the ring $S_j$ is free and
has a permissible system of parameters $(x,\frac{y}{x^j})$. In
particular,
$$
(X_1,Y_1)=\left(x,\frac{y}{x^{f_1}}\right)=
\left(\frac{x^{P_0}}{y^{P_{-1}}},\frac{y^{P_0}}{x^{P_1(f_1)}}\right)
$$
is a permissible system of regular parameters in $S_{F_1}$ with
$\n^*(X_1)=\m$ and $\n^*(Y_1)=\frac{r_2}{q}\m$. If $N=1$ then
$S=S_0\ra S_1\ra\dots\ra S_{F_N}$ is a sequence of free rings and
$S_{F_N}$ has a permissible system of parameters
$$
(X_N,Y_N)=(X_1,Y_1)=\left(x,\frac{y}{x^{f_1}}-c\right)=
\left(\frac{x^a}{y^b},\frac{y^q}{x^p}-c\right),
$$
where $c\in\k$ is the residue of $\frac{y^q}{x^p}$. Notice also that
$\n^*(X_N)=\m=(\n^*(x),\n^*(y))$.

If $N>2$ then for all $0<j\le f_2$, the ring $S_{F_1+j}$ is not free
and has a system of regular parameters $(\frac{X_1}{Y_1^j},Y_1)$. In
particular,
$$
(X_2,Y_2)=\left(\frac{X_1}{Y_1^{f_2}},Y_1\right)=
\left(\frac{x^{P_0+f_2P_1(f_1)}}{y^{f_2}},\frac{y^{P_0}}{x^{P_1(f_1)}}\right)=
\left(\frac{x^{P_2(f_1,f_2)}}{y^{P_1(f_2)}},
\frac{y^{P_0}}{x^{P_1(f_1)}}\right)
$$
are regular parameters in $S_{F_2}$ and $\n^*(X_2)=\frac{r_3}{q}\m$,
$\n^*(Y_2)=\frac{r_2}{q}\m$.

In general, for all $1<k<N$ and $0<j\le f_k$, the ring
$S_{F_{k-1}+j}$ is not free and has a system of regular parameters
$(\frac{X_{k-1}}{Y_{k-1}^j},Y_{k-1})$ if $k$ is even or
$(X_{k-1},\frac{Y_{k-1}}{X_{k-1}^j})$ if $k$ is odd. In particular,
if $k$ is even then $S_{F_k}$ has a system of regular parameters
$$
(X_k,Y_k)=\left(\frac{X_{k-1}}{Y_{k-1}^{f_k}},Y_{k-1}\right)
$$
where $\n^*(X_k)=\frac{r_{k+1}}{q}\m$ and
$\n^*(Y_k)=\frac{r_k}{q}\m$. We also notice that since
$$
\frac{X_{k-1}}{Y_{k-1}^{f_k}}=\frac{x^{P_{k-2}(f_1,\dots,f_{k-2})+f_kP_{k-1}(f_1,\dots,f_{k-1})}}
{y^{P_{k-3}(f_2,\dots,f_{k-2})+f_kP_{k-2}(f_2,\dots,f_{k-1})}}=
\frac{x^{P_k(f_1,\dots,f_k)}}{y^{P_{k-1}(f_2,\dots,f_k)}}
$$
the regular parameters $(X_k,Y_k)$ satisfy the equality
$$
(X_k,Y_k)=
\left(\frac{x^{P_k(f_1,\dots,f_k)}}{y^{P_{k-1}(f_2,\dots,f_k)}},
\frac{y^{P_{k-2}(f_2,\dots,f_{k-1})}}{x^{P_{k-1}(f_1,\dots,f_{k-1})}}\right).
$$

If $k$ is odd then $S_{F_k}$ has a system of regular parameters
$$
(X_k,Y_k)=\left(X_{k-1},\frac{Y_{k-1}}{X_{k-1}^{f_k}}\right)
$$
where $\n^*(X_k)=\frac{r_k}{q}\m$ and
$\n^*(Y_k)=\frac{r_{k+1}}{q}\m$. We notice that since
$$
\frac{Y_{k-1}}{X_{k-1}^{f_k}}=
\frac{y^{P_{k-3}(f_2,\dots,f_{k-2})+f_kP_{k-2}(f_2,\dots,f_{k-1})}}
{x^{P_{k-2}(f_1,\dots,f_{k-2})+f_kP_{k-1}(f_1,\dots,f_{k-1})}}=
\frac{y^{P_{k-1}(f_2,\dots,f_k)}}{x^{P_k(f_1,\dots,f_k)}}
$$
the regular parameters $(X_k,Y_k)$ satisfy the equality
$$
(X_k,Y_k)=
\left(\frac{x^{P_{k-1}(f_1,\dots,f_{k-1})}}{y^{P_{k-2}(f_2,\dots,f_{k-1})}},
\frac{y^{P_{k-1}(f_2,\dots,f_k)}}{x^{P_k(f_1,\dots,f_k)}}\right).
$$

Finally, if $N>1$ is odd then for all $0<j<f_N$ the ring
$S_{F_{N-1}+j}$ is not free and has a system of regular parameters
$(X_{N-1},\frac{Y_{N-1}}{X_{N-1}^j})$ . Moreover, $S_{F_N}$ is the
first free ring after a sequence of non-free rings
$S_{F_1+1}\ra\dots\ra S_{F_N-1}$. If $c\in\k$ is the residue of
$\frac{y^q}{x^p}$ then
\begin{multline*}
\quad(X_N,Y_N)=\left(X_{N-1},\frac{Y_{N-1}}{X_{N-1}^{f_N}}-c\right)=\\
=\left(\frac{x^{P_{N-1}(f_1,\dots,f_{N-1})}}{y^{P_{N-2}(f_2,\dots,f_{N-1})}},
\frac{y^{P_{N-1}(f_2,\dots,f_N)}}{x^{P_N(f_1,\dots,f_N)}}-c\right)=
\left(\frac{x^a}{y^b},\frac{y^q}{x^p}-c\right)\quad
\end{multline*}
form a permissible system of parameters in $S_{F_N}$ with
$\n^*(X_N)=\frac{r_N}{q}\m=\frac{1}{q}\m=(\n^*(x),\n^*(y))$.

If $N>1$ is even then for all $0<j<f_N$ the ring $S_{F_{N-1}+j}$ is
not free and has a system of regular parameters
$(\frac{X_{N-1}}{Y_{N-1}^j},Y_{N-1})$. Moreover, $S_{F_N}$ is the
first free ring after a sequence of non-free rings
$S_{F_1+1}\ra\dots\ra S_{F_N-1}$. If $c\in\k$ is the residue of
$\frac{y^q}{x^p}$ then
\begin{multline*}
(X_N,Y_N)=\left(\frac{X_{N-1}}{Y_{N-1}^{f_N-1}},
\frac{Y_{N-1}^{f_N}}{X_{N-1}}-c\right)=
\left(\frac{X_{N-1}}{Y_{N-1}^{f_N}}\cdot
Y_{N-1},\frac{Y_{N-1}^{f_N}}{X_{N-1}}-c\right)=\\
=\left(\frac{x^{P_N(f_1,\dots,f_N)-P_{N-1}(f_1,\dots,f_{N-1})}}
{y^{P_{N-1}(f_2,\dots,f_N)-P_{N-2}(f_2,\dots,f_{N-1})}},
\frac{y^{P_{N-1}(f_2,\dots,f_N)}}{x^{P_N(f_1,\dots,f_N)}}-c\right)
=\left(\frac{x^a}{y^b},\frac{y^q}{x^p}-c\right)
\end{multline*}
form a permissible system of parameters in $S_{F_N}$ with
$\n^*(X_N)=\frac{r_N}{q}\m=\frac{1}{q}\m=(\n^*(x),\n^*(y))$.

The following lemma summarizes the above discussion. We will often
refer to it in the rest of the paper.

\begin{lemma}\label{transformation}
Suppose that $S$ is a free ring and $(x,y)$ is a permissible system
of parameters in $S$ such that
$\dfrac{\n^*(y)}{\n^*(x)}=\dfrac{p}{q}$ for some coprime integers
$p$ and $q$. Let $k=\e(p,q)$,
$f_1=f_1(p,q)=\left[\dfrac{p}{q}\right]$ and let $a$ and $b$ be
nonnegative integers  such that $a\le p$, $b<q$, and $aq-bp=1$. Then
the sequence of quadratic transforms along $\n^*$
\begin{equation}\label{GSV}
S=S_0\ra S_1\ra\dots\ra S_{f_1}\ra S_{f_1+1}\ra\dots\ra S_{k-1}\ra
S_k
\end{equation}
has the following properties:
\begin{itemize}

\item[1)] $S_0,S_1,\dots,S_{f_1}$ and $S_k$ are free rings.

\item[2)] Non-free rings appear in {\rm (\ref{GSV})} if and only
if $k>f_1$,
 that is if $q>1$. In this case $S_{f_1+1},\dots,S_{k-1}$ are non-free.

\item[3)] $S_k$ has a permissible system of coordinates
$(X,Y)=\left(\dfrac{x^a}{y^b},\dfrac{y^q}{x^p}-c\right)$, where
$c\in\k$ is the residue of $\dfrac{y^q}{x^p}$. Moreover,
$\n^*(X)=(\n^*(x),\n^*(y))=\dfrac{1}{q}\n^*(x)$ and $x=X^q(Y+c)^b$,
$y=X^p(Y+c)^a$.
\end{itemize}
\end{lemma}
\begin{proof}
We only check that
$X^q(Y+c)^b=\dfrac{x^{aq}}{y^{bq}}\cdot\dfrac{y^{qb}}{x^{pb}}=x^{aq-bp}=x$
and
$X^p(Y+c)^a=\dfrac{x^{ap}}{y^{bp}}\cdot\dfrac{y^{qa}}{x^{pa}}=y^{aq-bp}=y$.
\end{proof}


\section{Properties of jumping polynomials}

In this section assumptions and notations are as in Section
\ref{jumpingpoly}. We fix regular parameters $(x,y)$ of $S$ and we
further assume that $(x,y)$ is a permissible system of parameters in
$S$ by setting $E$ to be the curve on $\spec S$ defined by $x=0$.



 For all $k>0$ let
$d_k=(p_1,p_2,\ldots,p_k).$ We will use this notation often in the
rest of the paper.

\begin{theorem}\label{strong}
Suppose that $R$ is a regular local ring dominated by $S$ and
$(u,v)$ are regular parameters of $R$ such that
\begin{align*}
u &=x^t\\
v &=y,
\end{align*}
where $t$ is a positive integer.

If $t|d_k$ for some $k>0$ then
$\{u,\{T_i\}_{i=1}^{k+1}\}=\{T'_i\}_{i=0}^{k+1}$ is the beginning of
a sequence of jumping polynomials in $R$. Moreover, for all $1\le
i\le k$ the pair of coprime integers defined in the construction of
jumping polynomials $\{T'_i\}_{i\ge 0}$ in $R$ is
$(p'_i,q'_i)=(\dfrac{p_i}{t},q_i)$.
\end{theorem}

\begin{proof}
Since $\n^*(u)=t$ in order to construct the sequence of jumping
polynomials $\{T'_i\}_{i\ge 0}$ in $R$ we use the following
valuation $\tilde\n$ of $K^*$:
$$
\tilde\n(f)=\frac{\n^*(f)}{t} \text { for all } f\in K^*.
$$
We have $T'_1=v=y=T_1$ and the coprime integers $p'_1$ and $q'_1$
are such that $\frac{p'_1}{q'_1}=\tilde\n(y)=\dfrac{p_1}{tq_1}$.
Assume $t|p_1$. Since $(p_1,q_1)=1$ we get $p'_1=\dfrac{p_1}{t}$ and
$q'_1=q_1$. Then $T'_2=v^{q_1}-\l'_1u^{p'_1}=y^{q_1}-\l'_1x^{p_1}$,
where $\l'_1$ is the residue of
$\dfrac{v^{q_1}}{u^{p'_1}}=\dfrac{y^{q_1}}{x^{p_1}}$, that is
$\l'_1=\l_1$ and $T'_2=T_2$. The statement is proved for $k=1$.

By induction on $k$ it suffices to show that the statement holds for
$k$ provided it holds for $k-1$. Then since $t|d_k$ and
$d_k|d_{k-1}$ by the inductive assumption we have $T'_i=T_i$ for all
$1\le i\le k$ and $(p'_i,q'_i)=(\dfrac{p_i}{t},q_i)$ for all $1\le
i\le k-1$. The coprime integers $p'_k$ and $q'_k$ satisfy the
following equality
$$
\frac{p'_k}{q'_k}=Q_{k-1}(\tilde\n(T_k)-q_{k-1}\tilde\n(T_{k-1}))=\frac{1}{t}
Q_{k-1}(\n^*(T_k)-q_{k-1}\n^*(T_{k-1}))=\frac{1}{t}\cdot\frac{p_k}{q_k}.
$$
Since $t|p_k$ and $(p_k,q_k)=1$ we get $p'_k=\dfrac{p_k}{t}$ and
$q'_k=q_k$. Then
$$
T'_{k+1}=(T'_k)^{q'_k}-\l'_k\prod_{i=0}^{k-1}(T'_i)^{n'_{k,i}}=
T_k^{q_k}-\l'_ku^{n'_{k,0}}\prod_{i=1}^{k-1}T_i^{n'_{k,i}},
$$
where $\l'_k$ is the residue of
$\dfrac{T_k^{q_k}}{u^{n'_{k,0}}\prod_{i=1}^{k-1}T_i^{n'_{k,i}}}$ and
$q_k\tilde\n(T_k)=n'_{k,0}+\sum_{i=1}^{k-1}n'_{k,i}\tilde\n(T_i)$
with $n'_{k,i}<q_i$ for all $1\le i\le k-1$. We notice that the last
equality is equivalent to
$$
q_k\n^*(T_k)=n'_{k,0}t+\sum_{i=1}^{k-1}n'_{k,i}\n^*(T_i).
$$
We also have
\begin{equation}\label{**}
q_k\n^*(T_k)=n_{k,0}+\sum_{i=1}^{k-1}n_{k,i}\n^*(T_i)
\end{equation}
from the construction of jumping polynomials in $S$. Thus from the
uniqueness of presentation (\ref{**}) we obtain $n'_{k,i}=n_{k,i}$
for all $1\le i\le k-1$ and $n'_{k,0}=\dfrac{n_{k,0}}{t}$. Thus,
$u^{n'_{k,0}}\prod_{i=1}^{k-1}T_i^{n'_{k,i}}=\prod_{i=0}^{k-1}T_i^{n_{k,i}}$,
$\l'_k=\l_k$, the residue of
$\dfrac{T_k^{q_k}}{\prod_{i=0}^{k-1}T_i^{n_{k,i}}}$, and
$T'_{k+1}=T_{k+1}$. This completes the proof.
\end{proof}

\begin{corollary}\label{xpowers}
For all $k>0$ all powers of $x$ that appear in
$T_2,T_3,\dots,T_{k+1}$ are multiples of $d_k$.
\end{corollary}

\begin{proof}
Notice that Theorem \ref{strong} in particular shows that if
$t|d_k$, then all powers of $x$ that appear in
$T_2,T_3,\dots,T_{k+1}$ are multiples of $t$.
\end{proof}

Our next goal is to describe the images of jumping polynomials under
blowups of $S$ along $\n^*$.

Before we state the next results we notice that if $S\subset\bar{S}$
is any subring of the $m_S$-adic completion $\hat{S}$ of $S$ we can
extend the valuation $\n^*$ to a valuation of $\bar{K}=QF(\bar{S})$
centered in $\bar{S}$. We first consider the unique extension
$\hat{\n}$ of $\n^*$ to $\hat{K}=QF(\hat{S})$ centered in $\hat{S}$,
then we restrict $\hat{\n}$ to $\bar{K}$. By abuse of notations we
will say that $\n^*$ is also a valuation of $\bar{K}$.

We will be mostly interested in the case where $\bar S$ is an
\'etale extension of $S$. If $n$ is a maximal ideal of $\bar S$, we
will say that the map $S \ra \bar S_{n}$ is {\it local \'etale}.
Most of the times we will have $n=m_{V^*} \cap \bar S$, the center
of the valuation.

We will first consider a sequence of ring extensions
$$
S=S'_0\ra \bar{S}'_0\ra S'_1\ra\dots\ra\bar{S}'_{i-1}\ra
S'_i\ra\dots
$$
such that for all $i\ge 0\;$ $\bar{S}'_i\ra S'_{i+1}$ is a quadratic
transform along $\n^*$ and $S'_i\ra\bar{S}'_i$ is a local \'etale
extension. As before let $E$ be a nonsingular irreducible curve on
$\spec S$, denote by $\p'_i$ the map $\spec S'_i\ra\spec S$ and by
$E'_i$ the reduced simple normal crossing divisor
${\p'}_i^{-1}(E)_{red}$.


In what follows, for all $i>0$ let $a_i,\, b_i$ be nonnegative
integers such that $a_iq_i-b_ip_i=1$ and $a_i\le p_i,\;b_i<q_i$. The
existence of $a_i$ and $b_i$ is due to the Euclidean division
algorithm. Let $k_0=0$ and $k_i = k_{i-1} + \e(p_i,q_i)$
($\e(p_i,q_i)$ was defined in Section \ref{arithmetic}).

\begin{lemma} \label{seq with etale}
There exists a sequence of ring extensions
\begin{multline*}
\quad S=S'_0\ra S'_1\ra\dots\ra S'_{k_1}\ra\bar{S}'_{k_1}\ra
S'_{k_1+1}\ra\dots\\
\dots\ra S'_{k_2}\ra\bar{S}'_{k_2}\ra\dots\ra
S'_{k_i}\ra\bar{S}'_{k_i}\ra\dots
\end{multline*}
such that for all $i > 0$ and $j\not=k_i$, $S'_j \rightarrow
S'_{j+1}$ and $\bar{S}'_{k_i}\ra S'_{k_i+1}$ are quadratic
transforms along $\n^*$, $\bar{S}'_{k_i}=S'_{k_i}[\a_i]_{m_{V^*}
\cap S'_{k_i}[\a_i]}$ are local \'etale extensions and the following
hold:
\begin{itemize}

\item[1)]$\a_i^{p_{i+1}}\in S'_{k_i}$ is a unit.

\item[2)]$S'_{k_i}$ is free and has a permissible system of
parameters $(z_i,w_i)$ such that $z_i$ is an exceptional parameter
and $w_i=\dfrac{T_{i+1}}{\prod_{j=0}^{i-1}T_j^{n_{i,j}}}$ is the
strict transform of $T_{i+1}$ in $S'_{k_i}$.

\item[3)] $\n^*(z_i)=\dfrac{1}{Q_i}\,$ and
$\n^*(w_i)=\dfrac{1}{Q_i}\cdot\dfrac{p_{i+1}}{q_{i+1}}$.

\item[4)] For all $0\le j\le i,\,$ $T_j=z_i^{Q_i\b_i}\t_{j,i}$,
where $\t_{j,i}\in S'_{k_i}$ is a unit.

\end{itemize}

\end{lemma}
\begin{proof}
We apply induction on $i$. For $i=1$, by Lemma \ref{transformation}
the ring $S'_{k_1}$ is free and has a system of regular parameters
$(z_1,w_1)$, where
$$
\begin{array}{ll}
z_1=\dfrac{x^{a_1}}{y^{b_1}} &
\quad\n^*(z_1)=(1,\dfrac{p_1}{q_1})=\dfrac{1}{q_1}\\
w_1=
\dfrac{y^{q_1}-\l_1x^{p_1}}{x^{p_1}}=\dfrac{T_2}{T_0^{n_{1,0}}}\quad
& \quad\n^*(w_1)=\b_2-p_1=\dfrac{1}{q_1}\cdot\dfrac{p_2}{q_2}.
\end{array}
$$
We also have $T_0=x=z_1^{q_1}(w_1+\l_1)^{b_1}=z_1^{q_1}\t_{0,1}$ and
$T_1=y=z_1^{p_1}(w_1+\l_1)^{a_1}=z_1^{p_1}\t_{1,1}$, where
$\t_{0,1}$ and $\t_{1,1}$ are units in $S'_{k_1}$.

Now assume that the lemma is true for $i-1$. We set
$\a_{i-1}=\left(\dfrac{\prod_{j=0}^{i-1}\t_{j,i-1}^{n_{i,j}}}
{\prod_{j=0}^{i-2}\t_{j,i-1}^{n_{i-1,j}q_i}}\right)^{\frac{1}{p_i}}$
and $\bar{S}'_{k_{i-1}}=S'_{k_{i-1}}[\a_{i-1}]_{m_{V^*} \cap
S'_{k_{i-1}}[\a_{i-1}]}.$ Then $\a_{i-1}^{p_i} \in S'_{k_{i-1}}$ and
$\bar{S}'_{k_{i-1}}$ is a local \'etale extension of $S'_{k_{i-1}}$.
Let $\bar{z}_{i-1}=z_{i-1}\a_{i-1}$, then $(\bar{z}_{i-1},w_{i-1})$
is a permissible system of parameters in $\bar{S}'_{k_{i-1}}$,
$\n^*(\bar{z}_{i-1})=\n(z_{i-1})=\dfrac{1}{Q_{i-1}}$ and
$\n^*(w_{i-1})=\dfrac{1}{Q_{i-1}}\cdot\dfrac{p_i}{q_i}$. Recall that
$k_i = k_{i-1} + \e(p_i,q_i)$. Therefore, by Lemma
\ref{transformation} the ring $S'_{k_i}$ is free and has a
permissible system of parameters $(z_i,w_i)$ such that
$$
\begin{array}{ll}
z_i=\dfrac{\bar{z}_{i-1}^{a_i}}{w_{i-1}^{b_i}}  \quad\text{ is an
exceptional parameter with }
\n^*(z_i)=\left(\dfrac{1}{Q_{i-1}},\dfrac{1}{Q_{i-1}}\cdot\dfrac{p_i}{q_i}
\right)=\dfrac{1}{Q_i},\\
w_i=\dfrac{w_{i-1}^{q_i}}{\bar{z}_{i-1}^{p_i}}-c_i, \text { where }
c_i\in\k \text{ is the residue of }
\dfrac{w_{i-1}^{q_i}}{\bar{z}_{i-1}^{p_i}}.
\end{array}
$$
We notice that for all $0\le l\le i-1$ the following equalities hold
\begin{align*}
\prod_{j=0}^l T_j^{n_{l+1,j}} & =\prod_{j=0}^l z_{i-1}^{Q_{i-1}\b_j
n_{l+1,j}}\t_{j,i-1}^{n_{l+1,j}}=z_{i-1}^{Q_{i-1}\sum_{j=0}^l
n_{l+1,j}\b_j}\prod_{j=0}^l\t_{j,i-1}^{n_{l+1,j}} \\
& =z_{i-1}^{Q_{i-1}q_{l+1}\b_{l+1}}
\prod_{j=0}^l\t_{j,i-1}^{n_{l+1,j}}=\bar{z}_{i-1}^{Q_{i-1}q_{l+1}\b_{l+1}}
\a_{i-1}^{-Q_{i-1}q_{l+1}\b_{l+1}}\prod_{j=0}^l\t_{j,i-1}^{n_{l+1,j}}.
\quad
\end{align*}
Hence,
$$
\prod_{j=0}^{i-1} T_j^{n_{i,j}}  = \bar{z}_{i-1}^{Q_i\b_i}
\a_{i-1}^{-Q_i\b_i}\prod_{j=0}^{i-1}\t_{j,i-1}^{n_{i,j}},\quad\quad
\prod_{j=0}^{i-2} T_j^{n_{i-1,j}} =
\bar{z}_{i-1}^{Q_{i-1}q_{i-1}\b_{i-1}}
\a_{i-1}^{-Q_{i-1}q_{i-1}\b_{i-1}}\prod_{j=0}^{i-2}\t_{j,i-1}^{n_{i-1,j}}
$$
and
\begin{align}
\dfrac{w_{i-1}^{q_i}}{\bar{z}_{i-1}^{p_i}} & = \dfrac{T_i^{q_i}}
{\bar{z}_{i-1}^{p_i}\prod_{j=0}^{i-2}T_j^{n_{i-1,j}q_i}}=
\dfrac{T_i^{q_i}}{\bar{z}_{i-1}^{Q_iq_{i-1}\b_{i-1}+p_i}
\a_{i-1}^{-Q_iq_{i-1}\b_{i-1}}\prod_{j=0}^{i-2}\t_{j,i-1}^{n_{i-1,j}q_i}}
\nonumber \\
& = \dfrac{T_i^{q_i}}{\bar{z}_{i-1}^{Q_i\b_i}
\a_{i-1}^{-Q_i\b_i+p_i}\prod_{j=0}^{i-2}\t_{j,i-1}^{n_{i-1,j}q_i}}=
\dfrac{T_i^{q_i}}{\prod_{j=0}^{i-1}
T_j^{n_{i,j}}}\cdot\dfrac{\prod_{j=0}^{i-1}\t_{j,i-1}^{n_{i,j}}}
{\prod_{j=0}^{i-2}\t_{j,i-1}^{n_{i-1,j}q_i}}
\a_{i-1}^{-p_i} \nonumber \\
& = \dfrac{T_i^{q_i}}{\prod_{j=0}^{i-1} T_j^{n_{i,j}}}. \label{find
wi}
\end{align}
Therefore, $c_i$ coincides with $\l_i$, and we get
$w_i=\dfrac{T_{i+1}}{\prod_{j=0}^{i-1} T_j^{n_{i,j}}}$ with
$\n^*(w_i)=\b_{i+1}-q_i\b_i=\dfrac{1}{Q_i}\cdot\dfrac{p_{i+1}}{q_{i+1}}$.

Finally, in view of Lemma \ref{transformation} we have
$z_{i-1}=z_i^{q_i}(w_i+\l_i)^{b_i}\a_{i-1}^{-1}$ and
$w_{i-1}=z_i^{p_i}(w_i+\l_i)^{a_i}$. Notice also that $Q_{i-1}\b_j$
is an integer for all $0\le j\le i-1$. Thus, for all $0\le j\le i-1$
we get
$$
T_j=z_{i-1}^{Q_{i-1}\b_j}\t_{j,i-1}=z_i^{Q_i\b_j}
((w_i+\l_i)^{b_i}\a_{i-1}^{-1})^{Q_{i-1}\b_j}\t_{j,i-1}=z_i^{Q_i\b_j}\t_{j,i},
$$
where $\t_{j,i}$ is a unit in $S'_{k_i}$, and
\begin{align*}
T_i & =
w_{i-1}\prod_{j=0}^{i-2}T_j^{n_{i-1,j}}=w_{i-1}z_{i-1}^{Q_{i-1}q_{i-1}\b_{i-1}}
\prod_{j=0}^{i-2}\t_{j,i-1}^{n_{i-1,j}}=\\
&= z_i^{Q_iq_{i-1}\b_{i-1}+p_i}(w_i+\l_i)^{a_i}
((w_i+\l_i)^{b_i}\a_{i-1}^{-1}
)^{Q_{i-1}q_{i-1}\b_{i-1}}\prod_{j=0}^{i-2}\t_{j,i-1}^{n_{i-1,j}}=z_i^{Q_i\b_i}\t_{i,i},
\end{align*}
where $\t_{i,i}$ is a unit in $S'_{k_i}$. This completes the proof
of the lemma.
\end{proof}

\begin{remark}\label{etale-lemma}
In our set-up, assume that $\bar S$ is a local \'etale extension of
$S$, $S'$ is a quadratic transform of $S$ along $\n^*$, and $\bar
S'$ is a quadratic transform of $\bar S$ along $\n^*$. Without loss
of generality, assume that $S$ has regular parameters $(x,y)$ with
$\n^*(y)\geq \n^*(x)$. Then $S'=S[x,\frac{y}{x}]_{(x,
\frac{y}{x}-\beta)}$, for some $\beta \in k$. Since $m_S\bar
S=m_{\bar S}$, we have that $\bar S$ has regular parameters $(x,y)$
and so $\bar S'=\bar S[x,\frac{y}{x}]_{(x, \frac{y}{x}-\beta')}$,
for some $\beta' \in k$. Since the quadratic transforms are along
$\n^*$, it follows that $\beta=\beta'$. Since $\bar S'$ is
essentially of finite type over $S'$ and $m_{S'}\bar S'=m_{\bar
S'}$, we have that $\bar S'$ is a local \'etale extension of $S'$.
Furthermore, if $\bar S=S[\alpha^{1/n}]_{m_{V^*}\cap
S[\alpha^{1/n}]}$, where $\alpha\in S$, we have that $\bar
S'=S'[\alpha^{1/n}]_{m_{V^*}\cap S'[\alpha^{1/n}]}$.

\end{remark}



\begin{theorem} \label{monoidalseq-lemma}
There exists a sequence of quadratic transforms along $\n^*$
$$
S=S_0\ra S_1\ra\dots\ra S_{k_1-1}\ra S_{k_1}\ra S_{k_1+1}\ra\dots\ra
S_{k_i}\ra\dots
$$
such that for all $i>0$, $S_{k_i}$ is free and has a system of
regular parameters $(x_i,y_i)$ such that $x_i$ is an exceptional
parameter, $T_j = x_i^{Q_i\b_j} \g_{j,i}$ for $0 \le j \le i$, where
$\g_{j,i} \in S_{k_i}$ is a unit, $\n^*(x_i)=\dfrac{1}{Q_i}\,$,
$y_i=\dfrac{T_{i+1}}{\prod_{j=0}^{i-1}T_j^{n_{i,j}}}$ is the strict
transform of $T_{i+1}$ in $S_{k_i}$ and
$\n^*(y_i)=\dfrac{1}{Q_i}\cdot\dfrac{p_{i+1}}{q_{i+1}}$.
\end{theorem}

\begin{proof}
We shall construct the required sequence from the sequence
$$
S = S_0'\ra S_1'\ra\dots\ra S_{k_1}'\ra\bar{S}_{k_1}'\ra S_{k_1+1}'
\ra\dots S_{k_2}'\ra\bar{S}_{k_2}'\ra\dots\ra
S_{k_i}'\ra\bar{S}_{k_i}'\ra\dots
$$
of Lemma \ref{seq with etale}. It suffices to show, by induction on
$i$, that we can construct a sequence
\begin{equation}
S = S_0 \ra S_1\ra\dots\ra S_{k_1}\ra S_{k_1+1}\ra\dots\ra S_{k_i}
\label{i-sequence}
\end{equation}
with the required properties and such that $\bar{S}'_{k_i} =
S_{k_i}[\a_1, \dots, \a_i]_{m_{V^*} \cap S_{k_i}[\a_1, \dots,
\a_i]}$ is a local \'etale extension of $S_{k_i}$, where
$\a_j^{p_{j+1}} \in S_{k_i}$ is a unit for $j=1, \dots, i$.

For $i = 1$, the sequence (\ref{i-sequence}) is given by taking $S_j
= S_j'$ for any $0 \le j \le k_1$ and setting $(x_1,y_1) =
(z_1,w_1)$. We also notice that $\bar{S}'_{k_1}=S_{k_1}[\a_1]_{
m_{V^*} \cap S_{k_1}[\a_1]}$, where $\a_1^{p_2}\in S_{k_1}$ is a
unit. In general, suppose that the sequence (\ref{i-sequence}) has
been constructed for $i-1$, i.e., we have a sequence of ring
extensions
\begin{align}
S=S_0 \rightarrow \dots \rightarrow S_{k_1} \rightarrow \dots
\rightarrow S_{k_{i-1}} \rightarrow \bar{S}_{k_{i-1}}' \rightarrow
S_{k_{i-1}+1}' \rightarrow \dots \rightarrow S_{k_i}' \rightarrow
\bar{S}_{k_i}' \rightarrow \dots \label{newseqi}
\end{align}
where $S_{j-1} \rightarrow S_j$ is a quadratic transform for $1 \le
j \le k_{i-1}$, $$\bar{S}_{k_{i-1}}' = S_{k_{i-1}}[\a_1, \dots,
\a_{i-1}]_{ m_{V^*} \cap S_{k_{i-1}}[\a_1, \dots, \a_{i-1}]}$$ is a
local \'etale extension (here, $\a_j^{p_{j+1}} \in S_{k_{i-1}}$ is a
unit for $j=1, \dots, i-1$) and $\bar{S}_{k_{i-1}}'\ra
S'_{k_{i-1}+1}$, $S'_{j-1}\ra S'_j$ is a quadratic transform for
$k_{i-1}+2\le j\le k_i$. By applying Remark \ref{etale-lemma} to the
subsequence
$$
\bar{S}_{k_{i-1}}'\ra S_{k_{i-1}+1}'\ra\dots\ra S_{k_i}'
$$
of the sequence (\ref{newseqi}) we obtain a new sequence of ring
extensions
$$
S=S_0 \rightarrow \dots \rightarrow S_{k_{i-1}} \rightarrow
S_{k_{i-1}+1} \rightarrow \dots \rightarrow S_{k_i} \rightarrow
\bar{S}_{k_i}' \rightarrow \dots
$$
where $S_{j-1} \rightarrow S_j$ is a quadratic transform for all $1
\le j \le k_{i}$ and $$\bar{S}'_{k_i} = S_{k_i}[\a_1, \dots, \a_i]_{
m_{V^*} \cap S_{k_i}[\a_1, \dots, \a_i]}$$ is a local \'etale
extension (here, $\a_j^{p_{j+1}} \in S_{k_{i}}$ is a unit for $j=1,
\dots, i$).

Let $(x_{i-1},y_{i-1})$ be a system of parameters in $S_{k_{i-1}}$
satisfying the required properties. Then by Lemma
\ref{transformation} we get that
$x_i=\dfrac{x_{i-1}^{a_i}}{y_{i-1}^{b_i}}$ is an exceptional
parameter of $S_{k_i}$ with $\n^*(x_i)=\dfrac{1}{Q_{i-1}}$ and
$\d=\dfrac{y_{i-1}^{q_i}}{x_{i-1}^{p_i}}$ is a unit in $S_{k_i}$.
Moreover, by computations similar to (\ref{find wi}), we have
$$
\dfrac{y_{i-1}^{q_i}}{x_{i-1}^{p_i}}
=\dfrac{T_i^{q_i}}{\prod_{j=0}^{i-1}
T_j^{n_{i,j}}}\cdot\dfrac{\prod_{j=0}^{i-1}\g_{j,i-1}^{n_{i,j}}}
{\prod_{j=0}^{i-2}\g_{j,i-1}^{n_{i-1,j}q_i}} =
\dfrac{T_i^{q_i}}{\prod_{j=0}^{i-1}T_j^{n_{i,j}}} \th,
$$
where $\th\in S_{k_i}$ is a unit.

This implies that
$\dfrac{T_i^{q_i}}{\prod_{j=0}^{i-1}T_j^{n_{i,j}}}\in S_{k_i}$ and,
therefore, $\dfrac{T_i^{q_i}}{\prod_{j=0}^{i-1}T_j^{n_{i,j}}}-\l_i=
\dfrac{T_{i+1}}{\prod_{j=0}^{i-1}T_j^{n_{i,j}}}=w_i\in S_{k_i}$.
Since $(z_i,w_i)$ form a permissible system of parameters in
$\bar{S}'_{k_i}$, by replacing the exceptional parameter $z_i$ by
$x_i$ we get a permissible system of parameters $(x_i,w_i)$ in
$S_{k_i}$. So we will choose a system of regular parameters
$(x_i,y_i)$ in $S_{k_i}$ by letting
$(x_i,y_i)=\left(\dfrac{x_{i-1}^{a_i}}{y_{i-1}^{b_i}},w_i\right)$.
Notice that $\n^*(y_i)=\dfrac{1}{Q_i}\cdot\dfrac{p_{i+1}}{q_{i+1}}$.

We further have $x_{i-1} = x_i^{q_i}\d^{b_i}$ and $y_{i-1} =
x_i^{p_i}\d^{a_i}$. Thus, recalling that $Q_{i-1}\b_j$ is an integer
for all $0 \le j \le i-1$, we obtain
$$
T_j = x_{i-1}^{Q_{i-1}\b_j}\g_{j,i-1} = x_i^{Q_i\b_j}
\d^{b_iQ_{i-1}\b_j}\g_{j,i-1}=x_i^{Q_i\b_j}\g_{j,i},
$$
where $\g_{j,i}\in S_{k_i}$ is a unit, and for $j = i$ we have
\begin{align*}
T_i & = y_{i-1}\prod_{j=0}^{i-2} T_j^{n_{i-1},j}= x_i^{p_i}\d^{a_i}
\left(\prod_{j=0}^{i-2} x_i^{n_{i-1,j}Q_i\b_j}
\g_{j,i}^{n_{i-1,j}}\right)=\\
& = x_i^{p_i + Q_i\sum_{j=0}^{i-2}
n_{i-1,j}\b_j}\left(\d^{a_i}\prod_{j=0}^{i-2}\g_{j,i}^{n_{i-1,j}}\right)=
x_i^{p_i + Q_iq_{i-1}\b_{i-1}}\g_{i,i}= x_i^{Q_i\b_i} \g_{i,i}
\end{align*}
where $\g_{i,i}\in S_{k_i}$ is a unit.

Hence, the sequence of ring extensions $S=S_0 \rightarrow S_1
\rightarrow \dots \rightarrow S_{k_i} \rightarrow \bar{S}_{k_i}'$
has required properties. The result is proved.
\end{proof}

\subsection{Remarks on generating sequences}\label{remgenseq}
 Assume that the value group of $\n^*$ is a
non-discrete subgroup of $\mathbb{Q}$. Suppose that
\begin{equation}\label{1}
S=S_0\ra S_1\ra S_2\ra\dots\ra S_j \ra\dots.
\end{equation}
is a sequence of quadratic transforms along $\n^*$ and $E_j$ is the
exceptional divisor on $S_j$ for all $j\ge 0$. Then a generating
sequence of $\n^*$ can be constructed as in \cite{C&P} (see also
\cite{Spi}, p. 150).

\begin{definition}
Set $s'_1=\bar{s}_0=0$. For all $i>0$ let $(s'_{i+1},\bar{s}_i)$ be
the pair of integers with the following properties:
\begin{itemize}
\item[1)] $\bar{s}_i$ is the biggest integer $s\ge s'_i$ such that
$S_{s'}$ is free for all $s'$ with $s'_i\le s'\le s$;

\item[2)] $s'_{i+1}$ is the smallest integer $s>\bar{s}_i$ such
that $S_s$ is free.
\end{itemize}
\end{definition}

We notice here that the set of free $S_j$ in (\ref{1}) is infinite,
as it follows from Theorem \ref{monoidalseq-lemma}. Thus the
sequences of integers $\{s'_i\}_{i>0}$ and $\{\bar{s}_i\}_{i\ge 0}$
are well defined.

\begin{remark}\label{olivierdef}
Let $\{Q_i\}_{i\ge 0}$ be a sequence of elements in $S$ such that
$Q_0$ is an exceptional parameter in $S$, $(Q_0,Q_1)$ form a system
of parameters in $S$ and the strict transform of div$(Q_1)$ in
$\spec S_{\bar{s}_1}$ is not empty. For each $i\ge 2$ let div$(Q_i)$
be an analytically irreducible curve in $\spec S$ such that the
strict transform of div$(Q_i)$ in $\spec S_{\bar{s}_i}$ is smooth
and transversal to $E_{\bar{s}_i}$. Then $\{Q_i\}_{i\ge 0}$ is a
generating sequence of $\n^*$ \cite{Spi}.
\end{remark}

We show that the set of all independent jumping polynomials in $S$
satisfies Remark \ref{olivierdef}. Therefore we have an alternative
argument that such set forms a generating sequence of $\n^*$. We
will need to use the irreducibility criterion of Cossart and
Moreno-Soc\'ias \cite[Theorem 6.2]{Cossart} in the form of Remark
7.17 of \cite{C&P}. Before we state this irreducibility criterion we
recall the notations of Section \ref{jumpingpoly}.

Suppose that $(x,y)$ are permissible parameters in $S$ and the value
group $\G^*$ of $\n^*$ is normalized so that $\n^*(x)=1$. Let
$\{T_i\}_{i\ge 0}$ denote the sequence of jumping polynomials in
$S$, $\{H_l\}_{l\ge 0}$ denote the sequence of independent jumping
polynomials in $S$ and $\{i_l\}_{l\ge 0}$ denote the sequence of
indexes such that $H_l=T_{i_l}$. Then for all $l\ge 1$ we have

\begin{equation}\label{Weierstrass}
H_{l+1}=H_l^{\bar{q}_l}-\l_{i_l}\prod_{j=0}^{l-1}H_j^{n_{i_l,i_j}}-
\sum_{i'=i_l+1}^{i_{l+1}-1}\l_{i'}\prod_{j=0}^l H_j^{n_{i',i_j}},
\end{equation}
where $0\le n_{i_l,i_j},n_{i',i_j}<\bar{q}_j$ for all
$i_l<i'<i_{l+1}$ and all $0\le j\le l$.

\begin{theorem}\label{irreducibility criterion}$($Remark 7.17,
\cite{C&P}$)$

Given a sequence of Weierstrass polynomials $\{H_l\}_{l\ge 1}$
satisfying (\ref{Weierstrass}) for all $l\ge 1$, set $\bar{\g}_0=1$
and define by induction on $l$ the values
$$
\bar{\g}_l=\dfrac{1}{\bar{q}_l}\sum_{j=0}^{l-1}n_{i_l,i_j}\bar{\g}_j.
$$
Let $\G_l=<\bar{\g}_0,\bar{\g}_1,\dots,\bar{\g}_l>$. Then
$\{H_l\}_{l\ge 0}$ is a generating sequence of a (uniquely
determined) valuation ring $\bar{V}$ of $\hat{S}=\k[[x,y]]$, whose
value group is a non-discrete subgroup of $\mathbb{Q}$ if for $l>0$
the $\bar{\g}_l$'s satisfy the following three properties:
\begin{itemize}

\item[1)] $\bar{q}_l\bar{\g}_l$ has order precisely $\bar{q}_l$ in
$\frac{\G_{l-1}}{\bar{q}_l\G_{l-1}}$,

\item[2)] $\bar{\g}_{l+1}>\bar{q}_l\bar{\g}_l$,

\item[3)] $\sum_{j=0}^{l}n_{i',i_j}\bar{\g}_j>\bar{q}_l\bar{\g}_l$
for all $i_l<i'<i_{l+1}$.
\end{itemize}
\end{theorem}

It follows from the construction of jumping polynomials that the
values $\bar{\g}_l$, defined for the sequence of independent jumping
polynomials, coincide with $\bar{\b}_l$ for all $l\ge 0$. Combining
Remark \ref{b*-inequality} and Corollary \ref{irreducibility
requirement} we see that independent jumping polynomials satisfy the
conditions of Theorem \ref{irreducibility criterion}. Thus the
sequence of independent jumping polynomials $\{H_l\}_{l\ge 0}$ is a
generating sequence for some valuation $\bar{\n}$ of $S$. In
particular, every element $H_l\in S$ is analytically irreducible in
$S$.

Let the sequence of quadratic transform of $S$
$$
S=S_0\ra S_1\ra S_2\ra\dots\ra S_{k_1}\ra\dots S_{k_2}\ra\dots\ra
S_{k_i}\ra\dots
$$
be as in Theorem \ref{monoidalseq-lemma}. Suppose that
$(x_{i-1},y_{i-1})$ are permissible regular parameters of
$S_{k_{i-1}}$, such that $\n^*(x_{i-1})=\frac{1}{Q_{i-1}}$ and
$\n^*(y_{i-1})=\frac{1}{Q_{i-1}}\cdot\frac{p_i}{q_i}$. Lemma
\ref{transformation} shows that non-free rings will appear in the
subsequence $S_{k_{i-1}}\ra\dots\ra S_{k_i}$ if and only if $q_i>1$,
that is if $\b_i$ is an independent j-value and $i=i_l$ for some
$l>0$. In this case the last free ring in this subsequence is
$S_{k_{i-1}+\left[\frac{p_i}{q_i}\right]}$ and $S_{k_i}$ is the
first free ring following
$S_{k_{i-1}+\left[\frac{p_i}{q_i}\right]}$. Thus for all $l>0$,
$\bar{s}_l=k_{i_l-1}+\left[\frac{p_{i_l}}{q_{i_l}}\right]$ and
$s'_{l+1}=k_{{i_l}}$.


Since $H_0=x$ and $H_1=y-\sum_{j=1}^{i_1-1}\l_jx^{\b_j}$,
$(H_0,H_1)$ form a permissible system of parameters in $S$.

By Theorem \ref{monoidalseq-lemma} there exists a permissible system
of parameters $(x_{i_l-1},y_{i_l-1})$ in $S_{k_{i_l-1}}$ such that
$y_{i_l-1}$ is the strict transform of $H_l$ in $S_{k_{i_l-1}}$. If
$f_1=\left[\frac{p_{i_l}}{q_{i_l}}\right]$ then
$S_{k_{i_l-1}+f_1}=S_{\bar{s}_l}$ and
$$
(X_{l},Y_{l})=(x_{i_l-1},\frac{{y_{i_l-1}}}{x_{i_l-1}^{f_1}})
$$
form a permissible system of parameters in $S_{\bar{s}_l}$. Thus the
strict transform of $H_l$ in $S_{\bar{s}_l}$ is $Y_l$. In
particular, the strict transform of div$(H_1)$ in $\spec
S_{\bar{s}_1}$ is not empty, and for $l\geq 2$, the strict transform
of div$(H_l)$ in $\spec S_{\bar{s}_l}$ is smooth and transversal to
$E_{\bar{s}_l}$. Therefore $\{H_l\}_{l\ge 0}\subset S$ form a
generating sequence of $\n^*$.

\begin{remark}
We have that $\{H_l\}_{l\ge 0}$ is a minimal generating sequence of
$\n^*$ if $p_1>1$. If $p_1=1$ then $\{H_l\}_{l>0}$ is a minimal
generating sequence of $\n^*$.
\end{remark}


\section{Monomialization of generating sequences}

The goal of this section is to prove the following theorem.

\begin{theorem}\label{monomialization}
Let $\k$ be an algebraically closed field of characteristic 0, and
let $K^*/K$ be a finite extension of algebraic function fields of
transcendence degree 2 over $\k$. Let $\n^*$ be a $\k$-valuation of
$K^*$, with valuation ring $V^*$ and value group $\G^*$, and let
$\n$ be the restriction of $\n^*$ to $K$, with valuation ring $V$
and value group $\G$. Suppose that $R \to S$ is an extension of
algebraic regular local rings with quotient fields $K$ and $K^*$
respectively, such that $V^*$ dominates $S$ and $S$ dominates $R$.
Then there exist sequences of quadratic transforms $R \to \bar{R}$
and $S \to \bar{S}$ along $\n^*$ such that $\bar{S}$ dominates
$\bar{R}$ and the map between generating sequences of $\n$ and
$\n^*$ in $\bar{R}$ and $\bar{S}$ respectively, has a toroidal
structure.
\end{theorem}

The lemma below is crucial in the proof of the theorem.

\begin{lemma}\label{chunk}
In the set up of Theorem \ref{monomialization}, assume that $\G^*$
is a subgroup of $\mathbb{Q}$ and $V^*/{m_{V^*}}=\k$. Suppose that
$R$ has regular parameters $(u,v)$ and $S$ has regular parameters
$(x,y)$ such that
\begin{equation}
\begin{array}{ll}
u &=  x^t \d \\
v &=  y,\\
\end{array}
\end{equation}
where $t$ is a positive integer and $\d$ is a unit in $S$. Let $p$
and $q$ be positive coprime integers such that
$\frac{\n^*(y)}{\n^*(x)}=\frac{p}{q}$ and let $k=\e(p,q)$. Let
$\bar{p}$ and $\bar{q}$ be positive coprime integers such that
$\frac{\n(v)}{\n(u)}=\frac{\bar{p}}{\bar{q}}$ and let
$\bar{k}=\e(\bar{p},\bar{q})$. Let $g$ be the greatest common
divisor of $t$ and $p$.

Then the sequences of quadratic transforms $R=R_0\ra R_1\ra\ldots\ra
R_{\bar{k}}$ and $S=S_0\ra S_1\ra\ldots\ra S_k$ along $\n^*$ satisfy
the following property: $R_{\bar{k}}$ and $S_k$ are free rings and
there exist permissible systems of regular parameters $(U,V)$ in
$R_{\bar{k}}$ and $(X,Y)$ in  $S_k$ such that
\begin{equation}
\begin{array}{ll}
U &=  X^{g} \D \\
V &=  Y \\
\end{array}
\end{equation}
for some unit $\D \in S_k$.
\end{lemma}

\begin{proof}
We notice first that
$\frac{\n(v)}{\n(u)}=\frac{\n^*(y)}{t\n^*(x)}=\frac{p}{tq}$. Writing
$t=gt'$ and  $p=gp'$, where $(t',p')=1$, gives $\bar p=p'$ and $\bar
q= qt'$. Also after possibly multiplying $u$ by a constant we may
assume that $\d=1+w$ for some $w\in m_S$.

Let $a$ and $b$ be nonnegative integers such that $a\leq p$, $b<q$
and $aq-bp=1$. Let $\bar a$ and $\bar b$ be nonnegative integers
such that $\bar a\leq \bar p$, $\bar b<\bar q$ and $\bar a \bar
q-\bar b\bar p=1$.

By Lemma \ref{transformation} applied to $S$ and $R$ respectively,
we get that $S_k$ has a permissible system of parameters
$(X,Y')=\left(\dfrac{x^a}{y^b},\dfrac{y^q}{x^p}-c\right)$, where
$c\in\k$ is the residue of $\dfrac{y^q}{x^p}$, and $R_{\bar{k}}$ has
a permissible system of parameters $(U,V)=\left(\dfrac{u^{\bar
a}}{v^{\bar b}},\dfrac{v^{\bar q}}{u^{\bar p}}-\bar c\right)$, where
$\bar c\in\k$ is the residue of $\dfrac{v^{\bar q}}{u^{\bar p}}$.
Moreover, $x=X^q(Y'+c)^b$, $y=X^p(Y'+c)^a$ and $u=U^{\bar q}(V+\bar
c)^{\bar b}$, $v=U^{\bar p}(V+\bar c)^{\bar a}$.

Now
$$
U=\frac{u^{\bar a}}{v^{\bar b}}=\frac{x^{t \bar a}\d^{\bar
a}}{y^{\bar b}}=\frac{[X^q(Y'+c)^b]^{t \bar a}\d^{\bar
a}}{[X^p(Y'+c)^a]^{\bar b}}=\frac{X^{qt \bar a}}{X^{p \bar
b}}\d^{\bar a}(Y'+c)^{bt\bar a- a\bar b}=X^g \D,
$$
where $\D=\d^{\bar a}(Y'+c)^{bt\bar a- a\bar b}$. Notice that the
last equality holds since $qt \bar a-p \bar b=qgt'\bar a-gp'\bar
b=g(qt'\bar a-p'\bar b)=g(\bar q \bar a-\bar p \bar b)=g$.

Furthermore notice that $\dfrac{v^{\bar q}}{u^{\bar
p}}=\dfrac{y^{\bar q}}{x^{t\bar p}}\d^{-\bar
p}=\dfrac{y^{qt'}}{x^{gt'p'}}\d^{-\bar p}=\left
(\dfrac{y^q}{x^p}\right)^{t'}\d^{-\bar p}=\left
(\dfrac{y^q}{x^p}\right)^{t'}(1+w)^{-\bar p}$. Therefore $\bar
c=c^{t'}$, and
$$
V=\frac {v^{\bar q}}{u^{\bar p}}-\bar c=\left
(\frac{y^q}{x^p}\right)^{t'}\d^{-\bar
p}-c^{t'}=(Y'+c)^{t'}(1+w)^{-\bar p}-c^{t'}=(Y'+c)^{t'}-c^{t'}+W,
$$
where $W\in wS\subset m_S$. Since $m_S\subset (X)S_k$ we have that
$V=Y'\d_2+XZ$ for some $Z\in S_k$ and unit $\d_2=\dfrac {\left
(\frac{y^q}{x^p}\right)^{t'}-c^{t'}}{\left
(\frac{y^q}{x^p}-c\right)}\in S_k$. Thus $(X,V)$ form a permissible
system of parameters in $S_k$. We set $Y=V$ to complete the proof.
\end{proof}

\begin{lemma}\label{main}
In the set up of Theorem \ref{monomialization}, assume that $\G^*$
is a subgroup of $\mathbb{Q}$ and $V^*/{m_{V^*}}=\k$. Suppose that
$R$ has regular parameters $(u,v)$ and $S$ has regular parameters
$(x,y)$ such that
\begin{equation}
\begin{array}{ll}
u &=  x^t \d \\
v &=  y,\\
\end{array}
\end{equation}
where $t$ is a positive integer and $\d$ is a unit in $S$.

Let $\tilde{S}=S[\d^{1/t}]_{ m_{V^*}\cap S[\d^{1/t}]}$ and let
$\{T_i\}_{i\ge 0}$ be the sequence of jumping polynomials in
$\tilde{S}$ such that $T_0=\tilde{x}=x\d^{1/t}$ and $T_1=y$. Then
either $\{u,\{T_i\}_{i>0}\}$ is a sequence of jumping polynomials in
$R$ or there exist sequences of quadratic transforms $R\ra R'$ and
$S\ra S'$ such that $R'$ has a system of regular parameters $(U,V)$,
$S'$ has a system of regular parameters $(X,Y)$ and
\begin{align*}
U &=X^g\D\\
V &=Y,
\end{align*}
where $g<t$ is a positive integer and $\D$ is a unit in $S'$.
\end{lemma}

\begin{proof} Without loss of generality we may assume that $\G^*$
is normalized so that $\n^*(x)=1$. Then $\n^*(\tilde{x})=1$ and
$\n(u)=t$.

For all $i>0$ let $p_i$ and $q_i$ be coprime integers defined in the
construction of jumping polynomials $\{T_i\}_{i\ge 0}$ in $\tilde
S$. Denote by $M=\min\{i>0|\;t\nmid p_i\}$. We assume first that
$M=\infty$, that is $p_i$ is multiple of $t$ for every $i$. Then
since $u=\tilde{x}^t$, by Theorem \ref{strong} we get that
$\{u,\{T_i\}_{i>0}\}$ is a sequence of jumping polynomials in $R$.

Assume now that $M<\infty$. Since $u=\tilde{x}^t$ and $t|p_i$ for
all $i<M$, by Theorem \ref{strong} we get that
$\{u,\{T_i\}_{i=1}^M\}=\{T'_i\}_{i=0}^M$ is the beginning of a
sequence of jumping polynomials in $R$ and for all $i< M$ the pairs
of coprime integers $(p'_i,q'_i)$ defined in the construction of the
sequence $\{T'_i\}_{i\ge 0}$ are $(p'_i,q'_i)=(\frac{p_i}{t},q_i)$.

Recall that the integers $k_i$ are defined as $k_0=0$ and
$k_{i}=k_{i-1}+\e(p_i,q_i)$ if $i>0$. Let $k'_0=0$ and
$k'_i=k'_{i-1}+\e(p'_i,q'_i)$ for all $i>0$. We will show first that
the sequences of quadratic transforms
$$
S=S_0\ra\dots\ra S_{k_1}\ra\dots\ra S_{k_{M-1}}
$$
and
$$
R=R_0\ra\dots\ra R_{k'_1}\ra\dots\ra R_{k'_{M-1}}
$$
have the following property: for all $0\le i\le M-1$ the rings
$R_{k'_i}$ and $S_{k_i}$ are free, there exist permissible systems
of parameters $(u_i,v_i)$ in $R_{k'_i}$ and $(x_i,y_i)$ in $S_{k_i}$
and a unit $\d_i\in S_{k_i}$ such that
$$
u_i=x_i^t\d_i,\quad\quad\quad v_i=y_i\quad
$$
and
$$
\quad\quad\n^*(x_i)=\frac{1}{Q_i},\quad\quad\quad
\n^*(y_i)=\frac{1}{Q_i}\cdot\frac{p_{i+1}}{q_{i+1}}.
$$

The statement is trivial for $i=0$. Assume that $i>0$ and that the
statement holds for $i-1$. Then Lemma \ref{chunk} applies to
$R_{k'_{i-1}}\subset S_{k_{i-1}}$. We notice that
$\dfrac{\n^*(y_{i-1})}{\n^*(x_{i-1})}=\dfrac{p_i}{q_i}$ and
$\dfrac{\n(v_{i-1})}{\n(u_{i-1})}=\dfrac{p_i}{tq_i}=\dfrac{p'_i}{q'_i}$,
and therefore $k=\e(p_i,q_i)$ and $\bar{k}=\e(p'_i,q'_i)$. Thus
$R_{k'_i}$ and $S_{k_i}$ are free rings and there exist permissible
systems of regular parameters $(u_i,w_i)$ in $R_{k'_i}$ and
$(x_i,z_i)$ in $S_{k_i}$ such that $u_i=x_i^t\d_i$ and $w_i=z_i$ for
some unit $\d_i\in S_{k_i}$.

Now by Theorem \ref{monoidalseq-lemma} applied to $R$ with $\n$
replaced by $\tilde{\n}=\frac{1}{t}\n$ we get that $R_{k'_i}$ has a
system of regular parameters $(h_i,v_i)$ such that $h_i$ is an
exceptional parameter,
$$
\n(h_i)=t\tilde{\n}(h_i)=t\dfrac{1}{q'_1\cdots
q'_i}=\dfrac{t}{q_1\cdots q_i}=\dfrac{t}{Q_i}
$$
and
$$
\n(v_i)=t\tilde{\n}(v_i)=t\dfrac{1}{q'_1\cdots
q'_i}\cdot\dfrac{p'_{i+1}}{q'_{i+1}}=
\dfrac{1}{Q_i}\cdot\dfrac{p_{i+1}}{q_{i+1}}.
$$
Since $u_i$ is also an exceptional parameter in $R_{k'_i}$ we have
$u_i=h_i\g$ for some unit $\g\in R_{k'_i}$. Therefore $(u_i,v_i)$
form a permissible system of parameters in $R_{k'_i}$ and
$\n(u_i)=\dfrac{t}{Q_i}$. Notice also that $v_i=\a u_i+\b w_i$,
where $\a,\b\in S_{k_i}$. Moreover, $\b$ is a unit in $R_{k'_i}$,
since the image of $v_i$ is a regular parameter in $R_{k'_i}/(u_i)$.
This implies that $v_i=\a x^t\d_i+\b z_i$ is also a regular
parameter in $S_{k_i}$ and $(x_i,v_i)$ form a permissible system of
parameters in $S_{k_i}$. We set $y_i=v_i$ and observe that
$\n^*(x_i)=\dfrac{1}{t}\n(u_i)=\dfrac{1}{Q_i}$ and
$\n^*(y_i)=\n(v_i)=\dfrac{1}{Q_i}\cdot\dfrac{p_{i+1}}{q_{i+1}}$.

To finish the proof of the lemma we apply Lemma \ref{chunk} to
$R_{k'_{M-1}}\subset S_{k_{M-1}}$. We have $p=p_M$, $q=q_M$ and
$\bar{p}=p'_M$, $\bar{q}=q'_M$. Thus $R'=R_{k'_M}$ has regular
parameters $(U,V)$ and $S'=S_{k_M}$ has regular parameters $(X,Y)$
such that
\begin{align*}
U &=X^g\D\\
V &=Y,
\end{align*}
where $\D$ is a unit in $S_{k_M}$ and $g=(p_M,t)<t$.
\end{proof}

We are now ready to prove Theorem \ref{monomialization}.

\begin{proof}
By the discussion of Section \ref{easycases} we only need to
consider the case when $\G^*$ is a subgroup of $\mathbb{Q}$ and
$\trdeg_{\k}(V^*/{m_{V^*}})=0$. Then by the Strong Monomialization
theorem we may assume that there exist regular parameters $(u,v)$ in
$R$ and $(x,y)$ in $S$ such that $u=x^t\d$ and $v=y$ for some unit
$\d\in S$. If $t=1$ then $R=S$ and the conclusion of the theorem is
trivial, so assume that $t>1$.

We set $\tilde{S}=S[\d^{1/t}]_{ m_{V^*}\cap S[\d^{1/t}]}$ and
$\tilde{x}=x\d^{1/t}$. Let $\{T_i\}_{i\ge 0}$ be a sequence of
jumping polynomials in $\tilde {S}$ such that $T_0=\tilde{x}$ and
$T_1=y$. For all $i>0$ let the coprime integers $p_i$ and $q_i$ be
defined as in the construction of jumping polynomials $\{T_i\}_{i\ge
0}$ in $\tilde{S}$.

First, let us assume that $t|p_i$ for all $i>0$. Then by the proof
of Lemma \ref{main} we have that $\{u,\{T_i\}_{i>0}\}$ is a sequence
of jumping polynomials in $R$. In particular, this implies that
$T_i\in S$ for all $i>0$.

If $\G^*$ is a discrete subgroup of $\mathbb{Q}$, after performing a
sequence of quadratic transforms along $\n^*$ and normalizing $\G^*$
we may assume that $\n^*(x)=1$ generates $\G^*$. In this case
Corollary \ref{groupvalue} shows that $q_i=1$ for all $i>0$. Then
$\n^*(T_i)=\sum_{j=1}^i p_j$ is a multiple of $t$ for all $i>0$.
Thus in view of Remark \ref{genvalues} we have that $\n^*(u)=t$
generates $\G$. By Theorem \ref{discretesequence} the sequence
$\{u,\{T_i\}_{i>0}\}$ form a generating sequence in $R$. Now since
$T_1=y$ and for all $i>0$
$$
T_{i+1}=v-\l_1u^{n_1}-\l_2u^{n_2}-\dots-\l_iu^{n_i}=
y-\l_1u^{n_1}-\l_2u^{n_2}-\dots-\l_iu^{n_i}
$$
is linear in $y$, repeating the proof of Theorem
\ref{discretesequence} we get that $\{x,\{T_i\}_{i>0}\}$ is a
generating sequence in $S$.

If $\G^*$ is a non-discrete subgroup of $\mathbb{Q}$, let
$\{H_i\}_{i\ge 0}$ be the sequence of independent jumping
polynomials in $\tilde{S}$. Then $\{u,\{H_i\}_{i>0}\}$ is a sequence
of independent jumping polynomials in $R$.  By Theorem
\ref{nondiscretesequence} we have that $\{H_i\}_{i\ge 0}$ is a
generating sequence in $\tilde{S}$ and $\{u,\{H_i\}_{i> 0}\}$ is a
generating sequence in $R$. Moreover, $\{H_i\}_{i\ge 0}$ is a
minimal generating sequence of $\n^*$ since $p_1>1$ as a multiple of
$t$. This implies that $\{x,\{H_i\}_{i>0}\}$ is a minimal generating
sequence in $S$.

If $t\nmid p_i$ for some $i>0$, let $R'$, $S'$ and $g$ be as in the
proof of Lemma \ref{main}. Notice that $g<t$, i.e., the exponent of
Strong Monomialization has dropped. We now repeat the above argument
starting with the rings $R'$ and $S'$ (instead of $R$ and $S$).
After a finite number of iterations we obtain the desired
conclusion.
\end{proof}

\begin{remark}
In the proof of Theorem \ref{monomialization} we have
$\{u,\{H_i\}_{i> 0}\}$ is a minimal generating sequence of $\n$ in
$R$ if $p_1\neq t$, otherwise $\{H_i\}_{i>0}$ is a minimal
generating sequence of $\n$ in $R$.
\end{remark}

\begin{remark}
Assumptions and notations are as in the statement of Theorem
\ref{monomialization}. By \cite[Theorem 6.1]{C&P} there exist
sequences of quadratic transforms $R \to \bar{R}$, $S \to \bar{S}$
along $\n^*$ such that $\bar R$ has regular parameters $(u,v)$,
$\bar S$ has regular parameters $(x,y)$ such that
\begin{equation}\label{ramification}
\begin{array}{ll}
u &=  x^e \d \\
v &=  y\\
\end{array}
\end{equation}

where $\d$ is a unit in $\bar S$, and $e=[\G^*:\G]$ is the
ramification index of $\n^*$ relative to $\n$. The monomial form of
(\ref{ramification}) is preserved by the sequences of quadratic
transforms of the proof of Theorem \ref{monomialization}.
Furthermore, the exponent $e$ does not drop under such sequences of
quadratic transforms (see the proof of \cite[Theorem 6.1]{C&P}). It
follows from the proof of Theorem \ref{monomialization} that the map
between generating sequences of $\n$ and $\n^*$ in $\bar{R}$ and
$\bar{S}$ respectively, has the desired toroidal structure.
\end{remark}



\begin{thebibliography}{10}

\bibitem{Ab} S. Abhyankar, {\it On the valuations centered in a domain}, Amer.
J. Math. 78 (1956), 321-348.

\bibitem{Cossart} V. Cossart, G. Moreno-Socias, {\it Racines approchees,
suites g\'{e}n\'{e}ratrices, suffisance des jets}, Valuation Theory
and its Applications II, F.-V. Kuhlmann, S. Kuhlmann, M. Marshall
editors, Fields Inst. Comm. 33, Amer. Math. Soc., Providence, RI,
361-459.

\bibitem{Dale2} S.D. Cutkosky, {\it Local factorization and
monomialization of morphisms}, Asterisque 260, 1999.

\bibitem{Dale3} S. D. Cutkosky, {\it Local factorization of birational
maps}, Adv. in Math. 132 (1997), 167-315.

\bibitem{C&P} S. D. Cutkosky, O. Piltant, {\it Ramification of Valuations},
Adv. in Math. 183 (2004), 1-79.

\bibitem{Dalebook} S.D. Cutkosky, Resolution of
Singularities, Graduate Studies in Mathematics 63, American
Mathematical Society, 2004.

\bibitem{FJ} C. Favre, M. Jonsson, {\it The valuative tree}, preprint, {\tt
arXiv:math.AC/0210265}.

\bibitem{GK} S. Greco, K. Kiyek, {\it General elements of complete ideals and
valuations centered at a two-dimensional regular local ring},
Algebra, Arithmetic and Geometry with Applications (West Lafayette,
IN, 2000), 381-455, Springer, Berlin, 2004.

\bibitem{MacL} S. MacLane, {\it A construction for absolute values
in polynomial rings}, Trans. Amer. Math. Soc. 40 (1936), no. 3,
363-395.

\bibitem{Spi} M. Spivakovsky, {\it Valuations in Functions Fields of
Surfaces}, Amer. J. Math. 112 (1990), 107-156.

\bibitem{T} B. Tessier, {\it Valuations, deformations and toric geometry},
Valuation Theory and its Applications II, F.-V. Kuhlmann, S.
Kuhlmann, M. Marshall editors, Fields Inst. Comm. 33, Amer. Math.
Soc., Providence, RI, 361-459.

\bibitem{Z} O. Zariski, {\it The reduction of the singularities of an
algebraic surface}, Annals of Math. 40 (1939), 639-689.

\bibitem{ZS} O. Zariski and P. Samuel, Commutative Algebra 2, Van
Nostrand, Princeton, 1960.
\end{thebibliography}
\end{document}